\documentclass[preprint,authoryear,12pt]{elsarticle}

\usepackage{amssymb,amsmath,amstext,psfrag}
\journal{International Journal of Plasticity}

\begin{document}

\begin{frontmatter}

%% Title, authors and addresses

%% use the tnoteref command within \title for footnotes;
%% use the tnotetext command for the associated footnote;
%% use the fnref command within \author or \address for footnotes;
%% use the fntext command for the associated footnote;
%% use the corref command within \author for corresponding author footnotes;
%% use the cortext command for the associated footnote;
%% use the ead command for the email address,
%% and the form \ead[url] for the home page:
%%
%% \title{Title\tnoteref{label1}}
%% \tnotetext[label1]{}
%% \author{Name\corref{cor1}\fnref{label2}}
%% \ead{email address}
%% \ead[url]{home page}
%% \fntext[label2]{}
%% \cortext[cor1]{}
%% \address{Address\fnref{label3}}
%% \fntext[label3]{}

\title{A viscoplasticity model with an enhanced control of the yield surface distortion}

%% use optional labels to link authors explicitly to addresses:
%% \author[label1,label2]{<author name>}
%% \address[label1]{<address>}
%% \address[label2]{<address>}

\author{A. V. Shutov, J. Ihlemann}

\address{Chemnitz University of Technology, Chair of Solid Mechanics, Chemnitz, Germany}

\begin{abstract}
%% Text of abstract
A new model of metal viscoplasticity, which takes combined
isotropic, kinematic, and distortional hardening into account, is
presented. The basic modeling assumptions are illustrated using a
new two-dimensional rheological analogy. This demonstrative
rheological model is used as a guideline for the construction of
constitutive relations. The nonlinear kinematic hardening is
captured using the well-known Armstrong-Frederick approach. The
distortion of the yield surface is described with the help of a
so-called distortional backstress. A distinctive feature of the
model is that any smooth convex saturated form of the yield surface
which is symmetric with respect to the loading direction can be
captured. In particular, an arbitrary sharpening of the saturated
yield locus in the loading direction combined with a flattening on
the opposite side can be covered. Moreover, the yield locus evolves
smoothly and its convexity is guaranteed at each hardening stage.
A strict proof of the thermodynamic consistency is provided.
Finally, the predictive capabilities of the material model are
verified using the experimental data for a very high work hardening
annealed aluminum alloy 1100 Al.

\end{abstract}

\begin{keyword}
viscoplasticity \sep yield function \sep kinematic hardening \sep
distortional hardening \sep rheology

%% keywords here, in the form: keyword \sep keyword

%% MSC codes here, in the form: \MSC code \sep code
%% or \MSC[2008] code \sep code (2000 is the default)
\MSC 74C10 \sep 74S05

\end{keyword}

\end{frontmatter}

% \linenumbers

%% main text
\section*{Nomenclature}
\begin{tabbing}
$\alpha$ \quad \quad \quad \quad \quad \quad    \= distortion parameter, cf. \eqref{DistortRheolo}, \eqref{DistortReal} \\
$\mathbf{X}_\text{k}$, $\mathbf{X}_\text{d}$ \>  backstress and distortional backstress, respectively \\
$R$ \>  isotropic hardening, cf. $\eqref{potent}_4$ \\
$\boldsymbol{\sigma}_{\text{eff}}$ \>  effective stress, cf. $\eqref{DefThet}_1$ \\
$\boldsymbol{R}_{\text{eff}}$ \> effective radial direction, cf. $\eqref{normFlow342}_2$ \\
$\boldsymbol{\varepsilon}_{\text{i}}$, $\boldsymbol{\varepsilon}_{\text{e}}$  \>
inelastic and elastic strains, respectively, cf. $\eqref{AddDec}_1$ \\
$\boldsymbol{\varepsilon}_{\text{ki}}$, $\boldsymbol{\varepsilon}_{\text{di}}$  \>
dissipative parts of $\boldsymbol{\varepsilon}_{\text{i}}$, cf. $\eqref{AddDec}_2$, $\eqref{AddDec}_3$ \\
$\boldsymbol{\varepsilon}_{\text{ke}}$, $\boldsymbol{\varepsilon}_{\text{de}}$   \>
conservative parts of $\boldsymbol{\varepsilon}_{\text{i}}$, cf. $\eqref{AddDec}_2$, $\eqref{AddDec}_3$ \\
$\lambda_{\text{i}}$ \>  inelastic multiplier, cf. $\eqref{Odqi}_2$ \\
$p$ \>  accumulated inelastic arc-length (Odqvist parameter), cf. $\eqref{Odqi}_3$ \\
$s$, $s_{\text{d}}$ \>  internal variables of isotropic hardening, cf. $\eqref{AddDec}_4$, \eqref{Odqvist} \\
$\theta$ \>  angle between $\boldsymbol{\sigma}^{\text{D}}_{\text{eff}}$
and $\mathbf{X}_\text{d}$, cf. $\eqref{DefThet}_2$ \\
$f$ \>  overstress, cf. \eqref{overDe}, \eqref{overDe3} \\
$\bar{f}$ \>  non-dimensional overstress, cf. \eqref{overdef} \\
$\text{El}(\bar{K}( \cdot , \alpha))$ \>  non-dimensional elastic domain, cf.  \eqref{DefEla} and Fig. \ref{fig3}b \\
$\bar{K}(\theta, \alpha)$ \>  non-dimensional yield stress function, cf. \eqref{DefK} \\
$\bar{K}^{\text{sat}}(\theta)$ \>  non-dimensional saturated yield stress function \\
$\mathcal{D}(\vec{y}, A)$ \> distance between $\vec{y} \in \mathbb{R}^2$ and $A \subset \mathbb{R}^2$, cf. \eqref{defdist}
\end{tabbing}

\section{Introduction}
\label{Intro}

Numerical simulation of complex metal forming operations is a
powerful tool to reduce the development costs and to optimize the
mechanical properties of a workpiece. For many
polycrystalline metals, the initial yield surface can be
approximated with sufficient accuracy by the conventional
Huber-Mises yield condition which implies initial plastic isotropy.
On the other hand, already very small plastic deformations may lead
to a significant change of the yield surface compared to the initial
state \citep{Annin2, Wegener, Dannemeyer, Steck,
KhanII, KhanIII}. It is well known that the residual stresses,
springback, damage evolution, and failure are highly dependent on the
accumulated plastic anisotropy of the material. In this work we
concentrate on the \emph{phenomenological} modeling of the plastic
anisotropy with especial emphasis on the distortional hardening.
A conventional approach to metal plasticity is used in the current study:
we suppose that a unique yield surface exists and that
the material behavior is purely elastic for
stresses within the yield surface.\footnote{Alternatively, different unconventional concepts
with numerous types of ``yield" surfaces exist. The so-called
``subloading surface models" (see \cite{Hashiguchi} and references therein)
allow to capture the plastic flow for stresses within the ``yield" surface.
Such approach allows to obtain a smooth transition from the elastic into the elasto-plastic range.
Moreover, some models with a smooth stress response can be constructed
using the concept of ``bounding surface"  (see, for instance, \cite{DafaliasP}).}

The state of the art phenomenological plasticity
is a result of accumulated efforts made by
generations of researchers. Unfortunately,
little academic credit was given to the paper written by \cite{Prager}.
Already in 1935, Prager combined the isotropic hardening of Odqvist type, the
distortional hardening for the prediction of cross hardening effect,
and the kinematic hardening for the Bauschinger effect.
Interestingly, the idea of modeling the Bauschinger effect by the kinematic translation
of the Huber-Mises yield surface in the stress space was taken by Prager
from a conference talk given by A. Reu\ss \ a year before - in 1934!

Within the classical phenomenological model of \cite{Chaboche1,
Chaboche2}, the isotropic expansion and kinematic translation of the
yield surface are considered, such that the yield surface is represented
by a hypersphere in the deviatoric stress
space.\footnote{Such yield surface can be represented by a
hypersphere in Ilyushin's space, as well (see \cite{Ilyushin}).}
Thus, the change of the form of the yield surface is neglected.
Such models can be used to simulate the stress response under
proportional loading. However, in general, the distortion of the yield
surface has to be considered under nonproportional loading with
abrupt change of the loading path.
Such loading conditions are typical not only for multi-stage forming processes,
but even for some single-stamping forming operations.
In order to control the rotation of a
hyperellipsoid which represents the yield surface withing a
Hill-type theory\footnote{The original approach of \cite{Hill} can be used to
describe a certain \emph{initial} plastic anisotropy, but not its evolution.},
\cite{Baltov} introduced a polynomial representation of the
corresponding 4th rank Hill-type anisotropy tensor in terms of the
strain tensor. According to \cite{Betten}, 4th and 6th rank
hardening tensors are postulated as functions of the plastic strain.
\cite{Dafalias} considered a general representation of the 4th rank tensor
as a polynomial function of the plastic strain.
In the paper by \cite{Helling}, the 4th rank Hill-type anisotropy
tensor is assumed to be a function of two backstress-like tensors.
In contrast to the above mentioned approaches, the approach of
Helling allows to take the dependence on the
strain path into account. For the same purpose, \cite{Rees,
Streilein, Kowalsky, Steck, DaScTs, FeigenbaumD, Noman, DafaliasFei, Pietryga}
and others modified the Chaboche-Rousselier model
introducing ordinary differential equations which describe the
evolution of tensor-valued internal variables of higher order (typically 4th
and 6th rank tensors). An alternative integral approach was
presented by \cite{Danilov}. In the paper of \cite{Grewolls},
evolution equations for higher order tensors were formulated in
integral form using the Danilov's approach. Both differential and
integral approaches mentioned above allow to take the dependence of
the hardening on the strain path. \cite{Kurtyka1, Kurtyka2} proposed
a geometric approach in order to simulate a complex distortion of
the yield surface.

The rigorous proof of convexity of the yield surface may become rather
difficult, if the 4th rank tensors are used \citep{Plesek}.
For instance, due to the complexity of the model presented by \cite{Pietryga},
the convexity of the yield surface was tested numerically by
plotting its two-dimensional projection at different
loading stages.

Probably, the most simple generalizations of the
Chaboche-Rousselier model are based on the use of second-rank
backstress-like tensors. Within this approach, the orientation of
the yield surface follows the loading path such that the change of
the loading direction leads to a reorientation of the yield surface
with a certain time lag. A short overview concerning different
approaches based on the use of backstress-like tensors is presented
by \cite{Wegener}. In particular, within the model of \cite{Ortiz},
the size of the elastic domain along a radial line emanating
from the origin of the yield surface
depends on the angle $\theta$ between
the effective deviatoric stress and the backstress direction. More precisely,
the critical norm of the effective deviatoric stress is given by a
Fourier cosine series of $\theta$. Thus, an arbitrary yield surface
which is symmetric with respect to the backstress direction can be
approximated. On the other hand, the convexity of the yield surface
imposes constraints on the Fourier coefficients. These constraints
complicate the construction of practical material models, especially
if the smooth evolution of the yield surface is intended. Another
special case was considered by \cite{Francois}. Within this
approach, certain egg-shaped yield surfaces can be modeled, such that
the egg-axis is oriented along a backstress-like tensor $\mathbf{X}_{\text{d}}$
and the degree of distortion is proportional to $\|\mathbf{X}_{\text{d}}\|$.
In particular, if Armstrong-Frederick type of hardening is used to describe
the evolution of $\mathbf{X}_{\text{d}}$, the distortion evolves in time
smoothly. The thermodynamic consistency was numerically tested by
Fran{\c c}ois. Next, within the approach presented by \cite{Panhans} as well as
\cite{Panhans2}, the distortional hardening was captured with the
help of a tensor-valued internal variable of the 2nd rank. The form
of the yield surface is given by the so-called lima{\c c}on  of
Pascal. Within the approach of Panhans it can be easily guaranteed
that the elastic domain is simply connected and convex. Later, in
the paper of \cite{ShutovPanh}, a two-dimensional rheological model
of distortional hardening was suggested, which implies the yield
surface to be the the lima{\c c}on of Pascal. This rheological model
was used to construct thermodynamically consistent constitutive
equations of finite strain plasticity/viscoplasticity.
In the paper by \cite{FeigenbaumD2}, an
existing material model was simplified in a thermodynamically
consistent manner such that second-rank backstress-like tensors are
used only.

%Crystal plasticity ????

A relatively new concept of
representative directions (see,
for example, \cite{Freund}) allows to generalize a uniaxial material model
to cover an arbitrary triaxial loading.
In order to compute the stress response, a numerical integration on the sphere
$S^2$ is required.
This concept, if combined with a uniaxial phenomenological
model of plasticity/viscoplasticity, can produce a new
phenomenological model with some realistic distortional hardening effects.
An interesting simplified approach to the description
of plastic anisotropy was proposed by \cite{Barlat}. Interestingly,
this approach does not include the concept of kinematic hardening explicitly, but some distortional
effects can be captured. The simplified approach to distortional hardening, which was developed by \cite{Aretz},
does not include kinematic hardening as well.
Further, we note that some models of crystal/polycrystal plasticity allow
the description of the yield surface distortion in a natural way (cf. \cite{Rousselier}).
For instance, in the paper by \cite{Fang}, the impact of
microstructural hardening parameters
on the form of the yield locus was analyzed in the finite strain context.
It was shown that for reduced latent hardening the yield surface exhibits
a larger curvature in the loading direction.

A new phenomenological model of metal plasticity is proposed in the
current study. The main features of the current model are as
follows:
\begin{itemize}
\item[\textbf{(i)}] a \textit{two-dimensional rheological
motivation} of constitutive equations, which provides insight into
main modeling assumptions;
\item[\textbf{(ii)}]  nonlinear isotropic
hardening of Voce type and nonlinear kinematic hardening
of Armstrong-Frederick type;
\item[\textbf{(iii)}] \emph{arbitrary} smooth
convex yield surface for saturated distortional hardening,
which is symmetric with respect to a backstress-like tensor $\mathbf{X}_{\text{d}}$;
\item[\textbf{(iv)}] degree of yield surface distortion is proportional to
$\|\mathbf{X}_{\text{d}}\|$; \emph{the convexity of the yield surface
is guaranteed} at each hardening stage;
\item[\textbf{(v)}] normality flow
rule; pressure-insensitive plasticity;
\item[\textbf{(vi)}] explicit formulation of the free energy density and
\emph{thermodynamic consistency};
\item[\textbf{(vii)}]
overstress type of viscoplasticity according to Perzyna rule.
\end{itemize}

In this paper, the temperature field is
assumed to be constant in time and
space.\footnote{The model is formulated in a thermodynamically
admissible manner. Therefore, its generalization to thermoplasticity
is straight-forward. The equation of heat conduction can be derived
directly from the energy balance, and an additional type of free
energy (so-called ``detached" free energy) can be introduced for
better prediction of temperature evolution, cf. \cite{ShutovIhle}.}
The model is formulated for infinitesimal strains
such that the extreme simplicity of the current approach is not
obscured by the geometric nonlinearities. At the same time, the
elegant technique of \cite{Lion}, which is based on the
consideration of rheological analogies can be used to generalize the
constitutive equations to finite strains
\citep{Helm, Shutov1, Henann, Vladimirov}. As it was shown by
\cite{ShutovPanh}, a similar technique can be implemented for
two-dimensional rheological models, as well.
Alternatively to the approach of Lion, the method of rheological models
proposed by \cite{Palmow} can be used to construct finite-strain
constitutive relations.

We conclude the introduction with a few remarks regarding notation.
The elements of $\mathbb{R}^2$ are denoted by
$\vec{x}$, $\vec{y}$. The notations
$\vec{x} \cdot \vec{y} := x_1 y_1 + x_2 y_2$ and
$ \|\vec{x}\| := \sqrt{\vec{x}
\cdot \vec{x}}$ stand for the scalar product and
the corresponding norm, respectively.
A coordinate-free tensor setting in $\mathbb{R}^3$
is implemented (cf. \cite{Itskov, Shutov2}).
Bold-faced symbols denote 1st- and 2nd-rank tensors in $\mathbb{R}^3$.
Superimposed dot denotes the material time derivative:
$\dot{x} = \frac{d}{d t} x$.
The symbol  `` $:$ " stands for the scalar product of two second-rank tensors
\begin{equation*}\label{scalprod}
\mathbf A : \mathbf B \ := \text{tr} (\mathbf A \cdot \mathbf B^{\text{T}}).
\end{equation*}
This scalar product gives rise to the
Frobenius norm as follows
\begin{equation*}\label{defi}
\| \mathbf A \| := \sqrt{ \mathbf A : \mathbf A}.
\end{equation*}
The identity tensor is denoted by $\mathbf 1$.
The notation $\mathbf A^{\text{D}}$ stands for a deviatoric part of a
tensor $\mathbf A^{\text{D}} := \mathbf A - \frac{1}{3}
\text{tr}(\mathbf A) \mathbf 1$.

\section{Rheological analogy}

\subsection{Two-dimensional rheological model}

Rheological models are useful for insight
into the aspects of material modeling.
Especially large body of information
is provided by rheological models if they are filled
with a physical content \citep{Petrov}.
Obviously, the conventional 1-dimensional rheological models are not suitable for
the description of the yield surface distortion.
Therefore, all considerations of this section are carried out in
two-dimensional space $\mathbb{R}^2$.
In the paper by \cite{ShutovPanh}, a two-dimensional rheological
model was suggested, which implies
that the yield surface is given by the lima{\c c}on  of
Pascal. \emph{A new extended rheological model}
of distortional hardening will be presented in this section.

\begin{figure}\centering
\psfrag{A}[m][][1][0]{$A$}
\psfrag{SV}[m][][1][0]{$m.StV$}
\psfrag{H1}[m][][1][0]{$H_{\text{ext}}$}
\psfrag{H2}[m][][1][0]{$H_{\text{dis}}$}
\psfrag{H3}[m][][1][0]{$H_{\text{kin}}$}
\psfrag{N1}[m][][1][0]{$m.N_{\text{dis}}$}
\psfrag{N2}[m][][1][0]{$m.N_{\text{kin}}$}
\psfrag{T}[m][][1][0]{$\theta$}
\psfrag{AA}[m][][1][0]{a)}
\psfrag{BB}[m][][1][0]{b)}
\scalebox{0.8}{\includegraphics{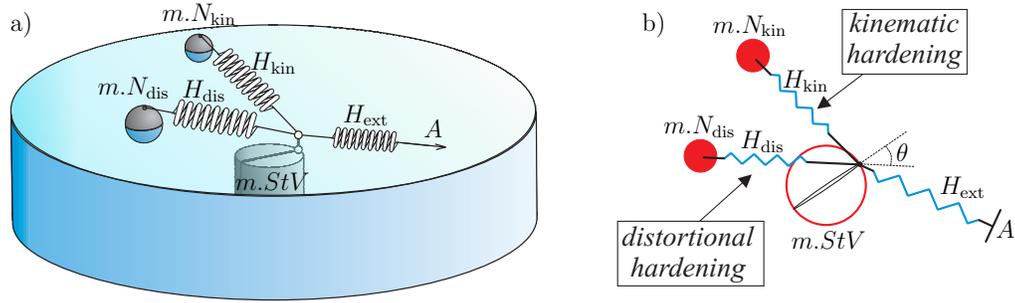}}
\caption{a) Two-dimensional rheological model
is built up of a modified St.-Venant element ($m.StV$),
Hooke-bodies $H_{\text{ext}}$, $H_{\text{kin}}$, $H_{\text{dis}}$, and
modified Newton elements $m.N_{\text{kin}}$ and $m.N_{\text{dis}}$;
b) Rheological model seen from above. The
angle between the ($m.StV$)-axis and $H_{\text{dis}}$ is
denoted by $\theta$.  \label{fig1}}
\end{figure}

\begin{figure}\centering
\psfrag{SV}[m][][1][0]{$m.StV$}
\psfrag{S}[m][][1][0]{$\vec{\sigma}_{H}$}
\psfrag{F}[m][][1][0]{$\vec{\sigma}_{N}$}
\psfrag{G}[m][][1][0]{$\vec{\sigma}$}
\psfrag{SM}[m][][1][0]{$-\vec{\sigma}_{H}$}
\psfrag{A}[m][][1][0]{$A$}
\psfrag{B}[m][][1][0]{$B$}
\psfrag{E}[m][][1][0]{$\vec{\varepsilon}_{H}$}
\psfrag{V}[m][][1][0]{$\frac{d}{d p}\vec{\varepsilon}_{N}$}
\psfrag{H}[m][][1][0]{$H$}
\psfrag{N}[m][][1][0]{$m.N$}
\psfrag{T}[m][][1][0]{$\theta$}
\psfrag{SE}[m][][1][0]{$\vec{\sigma}_{\text{eff}}$}
\psfrag{XK}[m][][1][0]{$-\vec{x}_{\text{k}}$}
\psfrag{XD}[m][][1][0]{$-\vec{x}_{\text{d}}$}
\psfrag{AA}[m][][1][0]{a)}
\psfrag{BB}[m][][1][0]{b)}
\psfrag{CC}[m][][1][0]{c)}
\scalebox{0.8}{\includegraphics{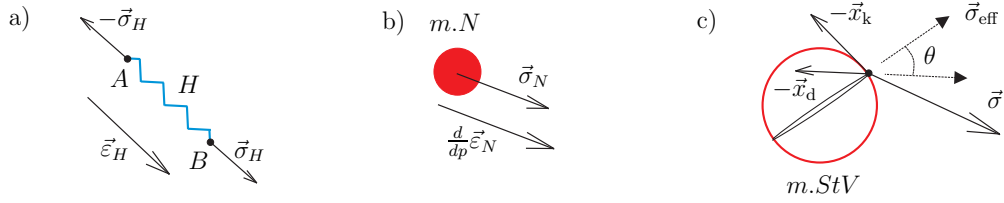}}
\caption{Behavior of idealized two-dimensional bodies: a) Hooke-body $H$. Its elongation is given by
$\vec{\varepsilon}_{H} = \protect\overrightarrow{AB} \in \mathbb{R}^2$ ;
b) Modified Newton-body $m.N$.
Inelastic arc-length $p$ is used instead of the physical time $t$ to formulate
the constitutive equations;
c) Modified St.-Venant element $m.StV$. The friction depends
on the angle $\theta$. \label{fig2}}
\end{figure}

In analogy to \cite{ShutovPanh}, we consider a mechanical system
which consists of a tank filled with a viscous fluid, a heavy solid
which rests on the flat bottom (modified St.-Venant element
$m.StV$), three elastic springs (Hooke-bodies $H_{\text{ext}}$,
$H_{\text{kin}}$, and $H_{\text{dis}}$) connected to the modified
St.-Venant element, and two spheres (modified Newton elements
$m.N_{\text{kin}}$ and $m.N_{\text{dis}}$) floating on the surface
of the fluid (Fig. \ref{fig1}).\footnote{An animated version of the
rheological model with only one modified Newton element is available
at http://www.youtube.com/watch?v=QEPc3pixbC0} The mechanical
properties of these idealized bodies are postulated as follows:
\begin{itemize}
\begin{item} ($H$): For the Hooke-bodies (see Fig. \ref{fig2}a), the spring force
$\vec{\sigma}_{H}$ is proportional to the length of the body,
and the force is oriented along the spring axis:
 $\vec{\sigma}_{H} = c \ \vec{\varepsilon}_{H}$.
Here, $\vec{\varepsilon}_{H} = \overrightarrow{AB} \in \mathbb{R}^2$,
and $c \geq 0$ is the stiffness of the spring.
In particular, the Hooke-body possesses a zero length in unloaded state.
\end{item}
\begin{item} ($m.N$): The two-dimensional Newton element is represented by
a sphere which is floating on the surface. Following the Newton's law of viscous flow, we assume that
the fluid resistance $\vec{\sigma}_{N}$ to the motion of the sphere is proportional to
its velocity $\frac{d}{d t}\vec{\varepsilon}_{N}$. Thus,
$\frac{d}{d t}\vec{\varepsilon}_{N} = \varkappa \ \vec{\sigma}_{N}$, where
$\varkappa \geq 0$ is a viscosity parameter.
Next, in order to obtain
the constitutive relations of the modified Newton-element,
the physical time $t$ is formally replaced
by the accumulated inelastic arc-length (Odqvist parameter) $p$.
Thus, we postulate for ($m.N$)-element
(see Fig. \ref{fig2}b):
\begin{equation}\label{NewtPos}
\frac{d}{d p}\vec{\varepsilon}_{N} = \varkappa \ \vec{\sigma}_{N}.
\end{equation}
Such modification is possible
whenever the inelastic arc-length $p$ is available.
The arc-length $p$ will be introduced formally
in the following. The use of this parameter
instead of the time $t$ allows to construct
rate-independent constitutive relations (see, for example, \cite{Haupt}).
\end{item}
\begin{item}
($m.StV$):
The heavy solid rests upon the bottom of the tank and there is a
friction between them. By $\vec{\sigma}$, $-\vec{x}_{\text{k}}$, and $-\vec{x}_{\text{d}}$
denote now the forces acting on the ($m.StV$)-element due to the elongation of
the Hooke-bodies $H_{\text{ext}}$, $H_{\text{kin}}$, and $H_{\text{dis}}$, respectively  (see Fig. \ref{fig2}c).
The force $\vec{\sigma}$ will be understood as an external load;
$\vec{x}_{\text{k}}$ and $\vec{x}_{\text{d}}$ will be responsible for the effects similar to
kinematic and distortional hardening, respectively.
The resulting (effective) force is thus given by
$\vec{\sigma}_{\text{eff}} = \vec{\sigma} - \vec{x}_{\text{k}} -
\vec{x}_{\text{d}}$.
Let the axis of the ($m.StV$)-element be always oriented along
the resulting (effective) force $\vec{\sigma}_{\text{eff}}$.
We suppose that the fluid resistance
opposed the rotation of the solid is negligible.
The ($m.StV$)-element remains at rest as long as $\| \vec{\sigma}_{\text{eff}} \| \leq \sqrt{2/3} K$, where
$\sqrt{2/3} K >0$ is a \emph{nonconstant} friction.
The function $K$ is computed as follows.
For $\vec{x}_{\text{d}}=\vec{0}$, we put
$K = K_0$, where $K_0$
is a given basic friction. Further, suppose
that $\vec{x}_{\text{d}} \neq \vec{0}$.
Let $\theta$ be the angle between the axis and $\vec{x}_{\text{d}}$:
$\theta = \text{arccos}\Big( \frac{\vec{\sigma}_{\text{eff}} \ \cdot \ \vec{x}_{\text{d}}}
{\|\vec{\sigma}_{\text{eff}}\| \ \|\vec{x}_{\text{d}}\|}\Big)$.
Moreover, let $\alpha = \|\vec{x}_{\text{d}}\| / x_d^{\text{max}}$ be a distortion parameter, which
is a unique function of $\|\vec{x}_{\text{d}}\|$.
Here, $x_d^{\text{max}} >0$ is the upper bound for $\|\vec{x}_{\text{d}}\|$,
therefore we get $\alpha \in [0, 1]$. Finally, we consider the friction to be a function of $\theta$ and $\alpha$:
$K = \bar{K}(\theta, \alpha) \ K_0$.
In particular, for a fixed $\vec{x}_{\text{d}} \neq \vec{0}$, the friction $K$ depends solely on the angle $\theta$.
A simple ansatz for $\bar{K}(\theta, \alpha)$ will be presented in the next subsection.
%such that $\bar{K}(\theta, \alpha) \rightarrow 1$ as $\alpha \rightarrow 0$.
\end{item}
\end{itemize}

\textbf{Remark 1.}
The choice of notations in this section is dictated by
the need to keep the structure of the rheological model
similar to the structure of small strain plasticity.
For that reason, the forces imposed on the
($m.StV$)-element by the Hooke-bodies
are denoted by $-\vec{x}_{\text{k}}$, $-\vec{x}_{\text{d}}$ rather than
$\vec{x}_{\text{k}}$, $\vec{x}_{\text{d}}$.

%Analogously to the one-dimensional approach,
%the body built up of a sphere $(N)$ connected in series to the elastic spring $(H)$
%is referred to as a Maxwell element.
%In the following the corresponding springs will be referred to as a Maxwell springs.

\subsection{Direction-dependent friction and definition of overstress}

\begin{figure}\centering
\psfrag{F}[m][][1][0]{$\alpha = 0$}
\psfrag{G}[m][][1][0]{$0 < \alpha < 1$}
\psfrag{H}[m][][1][0]{$\alpha = 1$}
\psfrag{EL}[m][][1][0]{$\text{El}(\bar{K}( \cdot , \alpha))$}
\psfrag{EO}[m][][1][0]{$\text{El}(\bar{K}( \cdot , 0))$}
\psfrag{U}[m][][1][0]{$\vec{e}_1$}
\psfrag{T}[m][][1][0]{$\theta$}
\psfrag{X}[m][][1][0]{$\vec{y}$}
\psfrag{x}[m][][1][0]{$\vec{y}$}
\psfrag{O}[m][][1][0]{$\vec{0}$}
\psfrag{ES}[m][][1][0]{$\text{El}^{\text{sat}}$}
%\psfrag{D}[m][][1][0]{$Distance(\vec{y}, A)$}
\psfrag{A}[m][][1][0]{$A$}
\psfrag{AA}[m][][1][0]{a)}
\psfrag{BB}[m][][1][0]{b)}
\psfrag{CC}[m][][1][0]{c)}
\scalebox{0.8}{\includegraphics{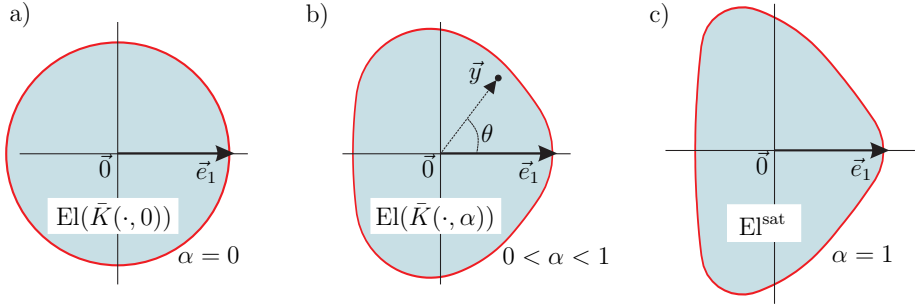}}
\caption{a) A closed unit disc which is used to represent the undistorted yield surface;
b) Geometric interpretation of $\bar{K}(\theta, \alpha)$ for a fixed $\alpha$.
For each $\alpha$ a closed set $\text{El}(\bar{K}( \cdot , \alpha)) \subset \mathbb{R}^2$ is defined;
c) A convex set $\text{El}^{\text{sat}}$ corresponds to the given function $\bar{K}^{\text{sat}}(\theta)$.
The boundary of $\text{El}^{\text{sat}}$ can be associated with a saturated form of the yield surface.\label{fig3}}
\end{figure}

Let $\vec{e}_1 = (1,0) \in \mathbb{R}^2$.
In this subsection we construct the non-dimensional function
$\bar{K}(\theta, \alpha)$ which plays a central role in the current study.
In terms of the rheological model introduced above,
this function is understood as a friction
coefficient, but in the following sections it
will be treated as a non-dimensional yield stress.
It is useful to interpret such functions geometrically
in terms of a parametric family of closed subsets in $\mathbb{R}^2$:
For each $\alpha \in [0,1]$ the corresponding subset $\text{El}(\bar{K}( \cdot , \alpha))$
consists of $\vec{y} \in \mathbb{R}^2$ such that $\| \vec{y} \| \leq  \bar{K}(\theta, \alpha)$, where
$\theta \in [0, \pi]$ is the angle between $\vec{y}$ and $\vec{e}_1$
(cf. Fig. \ref{fig3}b). More precisely, we put
\begin{equation}\label{DefEla}
\text{El}(\bar{K}( \cdot , \alpha)) := \big\{\vec{y} \in \mathbb{R}^2 / \{ \vec{0} \} \
: \ \| \vec{y} \| \leq  \bar{K}(\theta, \alpha), \ \text{where} \ \theta= \widehat{(\vec{y}, \vec{e}_1)}
 \big\} \cup \{ \vec{0} \},
\end{equation}
\begin{equation*}
\widehat{(\vec{y}, \vec{e}_1)} :=
\arccos\Big(\frac{\vec{y} \cdot \vec{e}_1}{ \|\vec{y}\| \ \|\vec{e}_1\| }\Big).
\end{equation*}
First, for $\alpha =0$, we postulate that $\bar{K}(\theta, 0) =1$, which implies that $\text{El}(\bar{K}( \cdot , 0))$ is
a closed unit disc in $\mathbb{R}^2$ (see Fig. \ref{fig3}a). Further, suppose that a smooth function
$\bar{K}^{\text{sat}}(\theta)$ is given such that
$\text{El}^{\text{sat}}:=\text{El}(\bar{K}^{\text{sat}}(\cdot))$ is convex and
$\bar{K}^{\text{sat}}(0)=1$ (see Fig. \ref{fig3}c). We need to find a smooth function
$\bar{K}(\theta, \alpha)$ such that $\text{El}(\bar{K}( \cdot , \alpha))$ will be convex
for all $\alpha$ and
\begin{equation*}
\bar{K}(0, \alpha) = 1 \ \text{for} \ \alpha \in [0,1], \quad \bar{K}(\theta, 1) = \bar{K}^{\text{sat}}(\theta) \
\text{for} \  \theta \in [0, \pi].
\end{equation*}

\textbf{Remark 2.}
As it will be shown in the following,
the given function $\bar{K}^{\text{sat}}(\theta)$ corresponds
to the form of a saturated distortional hardening with the maximum distortion (Fig. \ref{fig3}c).

\textbf{Remark 3.}
The parameter $\alpha$ should be understood as a distortion parameter, such that
$\alpha =0$ and $\alpha =1$ correspond to zero and maximum distortion, respectively.

In other words, an interpolation rule is needed between
the intact initial unit disc (corresponds to $\alpha=0$) and the maximum distorted set
$\text{El}^{\text{sat}}$ (corresponds to $\alpha=1$).

\textbf{Remark 4.}
Unfortunately, the linear interpolation rule
$\bar{K}(\theta, \alpha) = (1-\alpha) + \alpha \bar{K}^{\text{sat}}(\theta)$ is not suitable,
since, in general, the convexity of $\text{El}(\bar{K}( \cdot , \alpha))$ is violated
for some $\alpha \in (0,1)$.

\begin{figure}\centering
\psfrag{U}[m][][1][0]{$\vec{e}_1$}
\psfrag{EL}[m][][1][0]{$\text{El}(\bar{K}( \cdot , \alpha))$}
\psfrag{T}[m][][1][0]{$\theta$}
\psfrag{G}[m][][1][0]{$1 - \alpha$}
\psfrag{f}[m][][1][0]{$\bar{f}$}
\psfrag{x}[m][][1][0]{$\vec{y}$}
\psfrag{y}[m][][1][0]{$\vec{z}$}
\psfrag{R}[m][][1][0]{$\vec{y}(\theta, \alpha)$}
\psfrag{O}[m][][1][0]{$\vec{0}$}
\psfrag{AS}[m][][1][0]{$\alpha \text{El}^{\text{sat}}$}
\psfrag{D}[m][][1][0]{$\mathcal{D}(\vec{y}, A)$}
\psfrag{I}[m][][1][0]{$\{ \mathcal{D}(\vec{y}, \alpha \text{El}^{\text{sat}}) = 1 - \alpha \}$}
\psfrag{A}[m][][1][0]{$A$}
\psfrag{AA}[m][][1][0]{a)}
\psfrag{BB}[m][][1][0]{b)}
\psfrag{CC}[m][][1][0]{c)}
\scalebox{0.8}{\includegraphics{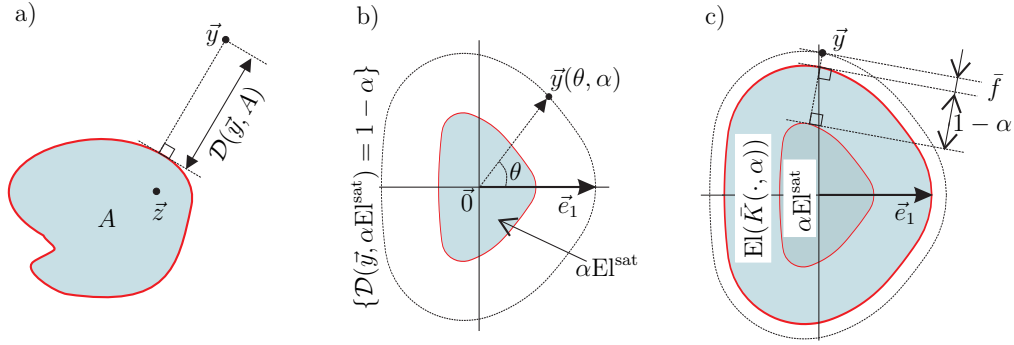}}
\caption{a) Definition of the distance $\mathcal{D}(\vec{y}, A)$ according to \eqref{defdist};
b) Definition of $\vec{y}(\theta, \alpha)$ for $\theta \in [0, \pi]$, $\alpha \in [0,1]$;
c) A line of constant overstress $\bar{f}$ for a fixed $\alpha \in [0,1]$.  \label{fig4}}
\end{figure}

The interpolation rule which is implemented in the current study
is constructed as follows.
First, for any $\vec{y} \in \mathbb{R}^2$, $A \subset \mathbb{R}^2$ we define the distance
in a natural way (see Fig. \ref{fig4}a)
\begin{equation}\label{defdist}
\mathcal{D}(\vec{y}, A) := \inf_{\vec{z} \in A}  \| \vec{y} - \vec{z} \|.
\end{equation}
Next, we define the product of $\alpha \in \mathbb{R}$ and $A \subset \mathbb{R}^2$ as
$\alpha A := \{ \alpha \vec{y}  \ : \ \vec{y} \in A \}$.
The set $\text{El}(\bar{K}( \cdot , \alpha))$ is obtained
from the set $\alpha \text{El}^{\text{sat}}$ by adding additional points whose distance
from $\alpha \text{El}^{\text{sat}}$ does not exceed $1-\alpha$ (see Fig. \ref{fig4}b)
\begin{equation}\label{interpSet}
\text{El}(\bar{K}( \cdot , \alpha)) := \big\{\vec{y} \in \mathbb{R}^2 :
\mathcal{D}(\vec{y}, \alpha \text{El}^{\text{sat}}) \leq 1-\alpha \big\}.
\end{equation}
Since $\text{El}^{\text{sat}}$ is convex, so is $\text{El}(\bar{K}( \cdot , \alpha))$.

Formally, since $\text{El}^{\text{sat}}$ is convex,
for each $\theta \in [0, \pi]$, and $\alpha \in [0,1]$
there exists a unique $\vec{y}(\theta, \alpha)$ such that
$\vec{y} = \| \vec{y} \| (\cos(\theta), \sin(\theta))$ and
$\mathcal{D}(\vec{y}, \alpha \text{El}^{\text{sat}}) = 1 - \alpha$  (see Fig. \ref{fig4}b).
Thus, in accordance with \eqref{interpSet}, we put
\begin{equation}\label{DefK}
\bar{K}(\theta, \alpha) := \| \vec{y}(\theta, \alpha) \|.
\end{equation}

In what follows, each set $\text{El}(\bar{K}( \cdot , \alpha))$
will be used to reflect the elastic region in the stress space.
Since a viscoplasticity model of overstress type is to be constructed,
a proper definition of the overstress will be needed.
For given $\vec{y} \in \mathbb{R}^2$, $\alpha \in [0,1]$ we define
a non-dimensional overstress as a distance from elastic domain:
\begin{equation}\label{overdef}
\bar{f}(\vec{y},\alpha) := \mathcal{D} \Big(\vec{y}, \text{El}(\bar{K}( \cdot , \alpha)) \Big) =
\big \langle \mathcal{D} \big(\vec{y}, \alpha \text{El}^{\text{sat}} \big) - (1-\alpha) \big \rangle,
\end{equation}
where $\langle x \rangle := \text{max}(x,0)$. The definition
is summarized in Fig. \ref{fig4}c.
Thus, for the numerical computation of the overstress it is sufficient to
evaluate $\mathcal{D} \big(\vec{y}, \alpha \text{El}^{\text{sat}} \big)$.
Such computation can be performed explicitly if the boundary of $\text{El}^{\text{sat}}$
is represented by a set of circular arcs.
A concrete algorithm is presented in Appendix A.
Note that such requirement is not restrictive,
since \emph{any smooth convex curve} can be approximated by circular arcs with sufficient accuracy.
In most cases 4 or 5 arcs are sufficient for practical purposes.

\textbf{Remark 5.}
The definition \eqref{overdef} of the overstress $\bar{f}$ will be advantageous
in connection with a normality flow rule. In particular,
the derivative of the overstress with respect to $\vec{y}$ possesses a unit norm for $\bar{f} > 0$, i.e.,
$ \| \frac{ \partial \bar{f}(\vec{y},\alpha)}{\partial \vec{y}} \| = 1$. Moreover, the set
$ \{ \vec{y} \ : \ \bar{f}(\vec{y}) \leq \bar{f}_0 \}$ is convex for all $\bar{f}_0 \geq 0$.

\begin{figure}\centering
\psfrag{U}[m][][1][0]{$\vec{e}_1$}
\psfrag{O}[m][][1][0]{$\vec{0}$}
\psfrag{ES}[m][][1][0]{$\text{El}^{\text{sat}}$}
\psfrag{A}[m][][1][0]{$\bar{K}^{\text{sat}}(\pi) <1  =   \bar{K}^{\text{sat}}(0)$}
\psfrag{B}[m][][1][0]{$\bar{K}^{\text{sat}}(\pi) = 1 = \bar{K}^{\text{sat}}(0)$}
\psfrag{C}[m][][1][0]{$\bar{K}^{\text{sat}}(\pi) > 1 = \bar{K}^{\text{sat}}(0)$}
\scalebox{0.8}{\includegraphics{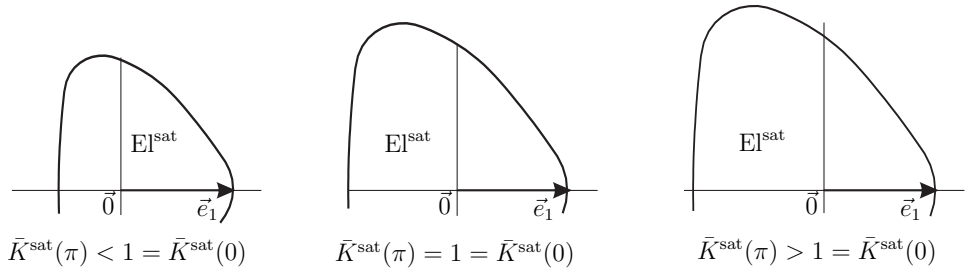}}
\caption{The material function $\bar{K}^{\text{sat}}(\theta)$ is not uniquely
determined by the form of the saturated yield surface.
The function is unique if $\bar{K}^{\text{sat}}(\pi) > 0$ is additionally specified.
\label{fig5}}
\end{figure}

\textbf{Remark 6.}
Using the interpolation rule proposed above,
the function $\bar{K}(\theta, \alpha)$ is uniquely defined by the material
function $\bar{K}^{\text{sat}}(\theta)$. Here, the
function $\bar{K}^{\text{sat}}(\theta)$ describes the saturated form of
the convex symmetric yield surface.
In some cases this form can be
identified experimentally (cf. Remark 9).
Figure \ref{fig5} demonstrates that the function $\bar{K}^{\text{sat}}(\theta)$ is not uniquely
determined even if the form of the saturated yield surface is known. This is due to the fact
the the position of the origin $\vec{0}$
relative to $\text{El}^{\text{sat}}$ is not unique.
The function $\bar{K}^{\text{sat}}(\theta)$ is uniquely determined
by specifying $\bar{K}^{\text{sat}}(\pi) > 0$, which
is a material parameter ``hidden" in the material function $\bar{K}^{\text{sat}}(\theta)$.
This parameter should be chosen in such a way that a better
description of the experimental data is achieved.
%For instance, if the origin $\vec{0}$ lies
%exactly between the current yield point ($\theta=0$) and the opposite yield point
%($\theta= \pi$), we have
%$\bar{K}^{\text{sat}}(\pi) = \bar{K}^{\text{sat}}(0) = 1$.
%In that case, the same stress response is predicted under proportional
%loading ($\theta=0$ or $\theta= \pi$) as by the model with kinematic hardening only.
%Alternatively, the position of the origin $\vec{0}$
%relatively to $\text{El}^{\text{sat}}$ can be chosen such that
%the area of $\text{El}^{\text{sat}}$ is equal to the area of the initial disc $\text{El}(\bar{K}( \cdot , 0))$.
%In that case, the area of the elastic domain is (almost) independent of $\alpha$.

\subsection{Some constitutive equations in two-dimensions}

Consider a system of (constitutive) equations as follows. The total
displacement of the point $A$ in Fig. \ref{fig1} with respect to the bottom will be denoted by $\vec{\varepsilon}$.
This displacement is a sum of the elastic elongation of the ($H_{\text{ext}}$)-body
and the inelastic displacement of the ($m.StV$)-body, denoted
by $\vec{\varepsilon}_{\text{e}}$ and $\vec{\varepsilon}_{\text{i}}$, respectively:
\begin{equation}\label{BDecomReol}
\vec{\varepsilon} = \vec{\varepsilon}_{\text{e}} +
\vec{\varepsilon}_{\text{i}}.
\end{equation}
The displacement of the ($m.StV$)-body, in turn, is composed of the elastic ($H_{\text{kin}}$)-elongation
and the inelastic ($m.N_{\text{kin}}$)-displacement. Analogous decomposition
holds for $H_{\text{dis}}$ and $m.N_{\text{dis}}$, as well:
\begin{equation}\label{DecomReol}
\vec{\varepsilon}_{\text{i}} = \vec{\varepsilon}_{\text{ke}} + \vec{\varepsilon}_{\text{ki}},
\quad
\vec{\varepsilon}_{\text{i}} = \vec{\varepsilon}_{\text{de}} + \vec{\varepsilon}_{\text{di}}.
\end{equation}
The total potential energy of the system equals
\begin{equation*}
\psi = \psi_{\text{el}} (\vec{\varepsilon}_{\text{e}}) +
\psi_{\text{kin}} (\vec{\varepsilon}_{\text{ke}}) +
\psi_{\text{dis}} (\vec{\varepsilon}_{\text{de}}) = \mu \ \| \vec{\varepsilon}_{\text{e}} \|^2+ \frac{c_{\text{k}}}{2}
 \ \|\vec{\varepsilon}_{\text{ke}} \|^2 + \frac{c_{\text{d}}}{2} \ \|\vec{\varepsilon}_{\text{de}} \|^2.
\end{equation*}
Here, $\mu , c_{\text{k}},  c_{\text{d}} \geq 0$ are the stiffnesses of $H_{\text{ext}}$,
$H_{\text{kin}}$, and $H_{\text{dis}}$, respectively. For the forces $\vec{\sigma}$, $\vec{x}_{\text{k}}$, and $\vec{x}_{\text{d}}$ we get
\begin{equation}\label{stresses}
\vec{\sigma} = \frac{\partial \psi_{\text{el}} (\vec{\varepsilon}_{\text{e}})}{\partial \vec{\varepsilon}_{\text{e}}}, \quad
\vec{x}_{\text{k}} = \frac{\partial \psi_{\text{kin}} (\vec{\varepsilon}_{\text{ke}})}{\partial \vec{\varepsilon}_{\text{ke}}}, \quad
\vec{x}_{\text{d}} = \frac{\partial \psi_{\text{dis}} (\vec{\varepsilon}_{\text{de}})}{\partial \vec{\varepsilon}_{\text{de}}}.
\end{equation}
The overstress $f$ is defined as a function of $\| \vec{\sigma}_{\text{eff}} \| =  \| \vec{\sigma} - \vec{x}_{\text{k}} -
\vec{x}_{\text{d}} \|$, $\theta = \text{arccos}\Big( \frac{\vec{\sigma}_{\text{eff}} \ \cdot \ \vec{x}_{\text{d}}}
{\|\vec{\sigma}_{\text{eff}}\| \ \|\vec{x}_{\text{d}}\|}\Big)$, and the
distortion parameter $\alpha = \|\vec{x}_{\text{d}}\| / x_d^{\text{max}}$ by
\begin{equation}\label{overDe}
f(\vec{\sigma}, \vec{x}_{\text{k}}, \vec{x}_{\text{d}})
=\tilde{f}(\| \vec{\sigma}_{\text{eff}} \|, \theta, \alpha) := \sqrt{\frac{2}{3}} K_0 \ \bar{f}(\vec{y},\alpha),
\end{equation}
\begin{equation}\label{overDe2}
\vec{y} :=\frac{ \| \vec{\sigma}_{\text{eff}} \|}{\sqrt{2/3} K_0} (\cos(\theta), \sin(\theta)).
\end{equation}
Due to the fact that $\| \frac{ \partial \bar{f}(\vec{y},\alpha)}{\partial \vec{y}} \| = 1$ for positive overstress, we have
\begin{equation*}
\Big\| \frac{ \partial f(\vec{\sigma}, \vec{x}_{\text{k}}, \vec{x}_{\text{d}}) }{\partial \vec{\sigma}} \Big\| = 1, \ \text{for} \
f>0.
\end{equation*}
We postulate the normality flow rule (normality to the hypersurface of constant overstress) in
combination with the Perzyna-type of viscoplasticity \citep{Perzyna2}
\begin{equation*}
\dot{\vec{\varepsilon}}_{\text{i}} = \lambda_{\text{i}} \frac{ \partial f(\vec{\sigma}, \vec{x}_{\text{k}},
\vec{x}_{\text{d}}) }{\partial \vec{\sigma}},
\ \text{for} \ f>0, \quad
\dot{\vec{\varepsilon}}_{\text{i}} = \vec{0} \ \text{for} \ f = 0; \quad
\lambda_{\text{i}} = \frac{\displaystyle 1}{\displaystyle
\eta}\Big( \frac{\displaystyle 1}{\displaystyle k_0} f \Big)^{m}.
\end{equation*}
Here, $\eta > 0$ and $m \geq 1$ are parameters of the Perzyna rule;
$k_0 >0$ is used to get a dimensionless term in the parentheses.

Note that the forces $\vec{x}_{\text{k}}$ and $\vec{x}_{\text{d}}$
act on the modified Newton-elements $m.N_{\text{kin}}$ and $m.N_{\text{dis}}$,
respectively. Thus, we get in accordance with \eqref{NewtPos}
\begin{equation*}\label{saturIni}
\frac{d \ \vec{\varepsilon}_{\text{ki}}}{d  p} =  \ \varkappa_\text{k} \ \vec{x}_{\text{k}}, \quad
\frac{d \ \vec{\varepsilon}_{\text{di}}}{d  p} =  \ \varkappa_\text{d} \ \vec{x}_{\text{d}}.
\end{equation*}
Here, $\varkappa_\text{k}, \varkappa_\text{d} \geq 0$ are modified viscosity parameters describing
$m.N_{\text{kin}}$ and $m.N_{\text{dis}}$, respectively.
 Let the evolution of the inelastic arc-length $p$ be given by
$\dot{p} = \lambda_{\text{i}} = \| \dot{\vec{\varepsilon}}_{\text{i}} \|$. Thus, we obtain
\begin{equation}\label{satur}
\dot{\vec{\varepsilon}}_{\text{ki}} = \lambda_{\text{i}} \ \varkappa_\text{k} \ \vec{x}_{\text{k}}, \quad
\dot{\vec{\varepsilon}}_{\text{di}} = \lambda_{\text{i}} \ \varkappa_\text{d} \ \vec{x}_{\text{d}}.
\end{equation}
It follows from $\eqref{stresses}_3$ and $\eqref{satur}_2$ that for proper initial
conditions we have
$\|\vec{x}_{\text{d}}\| \leq 1/\varkappa_{\text{d}}$.
By putting $ x_d^{\text{max}} = 1/\varkappa_{\text{d}}$ we specify
the definition of the distortion parameter $\alpha$ (cf. Section 2.1)
\begin{equation}\label{DistortRheolo}
\alpha := \varkappa_{\text{d}} \ \|\vec{x}_{\text{d}}\|.
\end{equation}
Equations \eqref{satur} in combination
with $\eqref{stresses}_2$ and $\eqref{stresses}_3$ describe the evolution of the ``backstresses" in the
hardening/recovery format. The saturation of $\vec{x}_{\text{d}}$
implies the saturation of the ``distortional hardening", which takes
place much faster than the saturation of the ``kinematic hardening".
Thus,  ``slow" and  ``fast" saturation should be assumed for
$\vec{x}_{\text{k}}$ and $\vec{x}_{\text{d}}$, respectively.
Note that the dependence on the strain path is captured in a vivid way, such that
the system exhibits fading memory: only
the most recent part of the $\vec{\varepsilon}$-path influences the current  ``stress" state $\vec{\sigma}$.

\section{Material model of viscoplasticity}

\subsection{Closed system of constitutive equations}

Let us formulate a system of constitutive equations of viscoplasticity.
First, we suppose that the volumetric response is elastic. More precisely,
the hydrostatic stress component
$\text{tr}  \boldsymbol{\sigma}$ is assumed to be a linear function
of $\text{tr}  \boldsymbol{\varepsilon}$.
Next, suppose that the deviatoric stress component $\boldsymbol{\sigma}^{\text{D}}$ depends solely on the history
of the strain deviator $\boldsymbol{\varepsilon}^{\text{D}}$. In order to
describe this dependence, we generalize
the two-dimensional constitutive equations presented in the previous
section to five dimensions.\footnote{Mathematically,
$\boldsymbol{\sigma}^{\text{D}}$ and $\boldsymbol{\varepsilon}^{\text{D}}$ are
elements of a 5-dimensional vector space
of trace-free (deviatoric) second rank symmetric tensors.}
During the generalization we have to make sure that the resulting
model inherits the properties of the two-dimensional rheological model.
Toward that end, the displacements
and forces are formally replaced by
deviatoric strains and stresses, respectively;
the scalar product in $\mathbb{R}^2$ is
replaced by the scalar product of two second-rank
tensors.\footnote{Using a similar approach, a two-dimensional rheological
model was already generalized
to cover \emph{finite strain} vsicoplasticity \citep{ShutovPanh}.}
In order to take the isotropic hardening into account,
the constant parameter $K_0$ is now formally replaced by
$K_0 + R$, where $R$ is a hardening variable. In order to describe the evolution of $R$,
we introduce a scalar strain-like internal variable $s$ (which is similar to the inelastic arc-length $p$),
its dissipative part $s_\text{d}$, and its conservative part $s_\text{e}$.

For the strain tensor
$\boldsymbol{\varepsilon}$ consider its inelastic part $\boldsymbol{\varepsilon}_{\text{i}}$ and elastic part
$\boldsymbol{\varepsilon}_{\text{e}}$.
Let $\boldsymbol{\varepsilon}_{\text{ki}}$ and $\boldsymbol{\varepsilon}_{\text{ke}}$
be the dissipative and conservative parts of $\boldsymbol{\varepsilon}_{\text{i}}$, which are
connected to the nonlinear kinematic hardening. Analogously,
$\boldsymbol{\varepsilon}_{\text{di}}$ and $\boldsymbol{\varepsilon}_{\text{de}}$ are parts of
$\boldsymbol{\varepsilon}_{\text{i}}$ associated to the distortional hardening. More precisely, we postulate
\begin{equation}\label{AddDec}
\boldsymbol{\varepsilon} = \boldsymbol{\varepsilon}_{\text{e}} + \boldsymbol{\varepsilon}_{\text{i}}, \quad
\boldsymbol{\varepsilon}_{\text{i}} = \boldsymbol{\varepsilon}_{\text{ke}} + \boldsymbol{\varepsilon}_{\text{ki}}, \quad
\boldsymbol{\varepsilon}_{\text{i}} = \boldsymbol{\varepsilon}_{\text{de}} + \boldsymbol{\varepsilon}_{\text{di}}, \quad
s = s_\text{e} + s_\text{d}.
\end{equation}
Note that the first decomposition is related to \eqref{BDecomReol}.
Moreover, $\eqref{AddDec}_2$ and $\eqref{AddDec}_3$ can be motivated
by \eqref{DecomReol}. The evolution of the state
of the material is captured by
the inelastic flow $\dot{\boldsymbol{\varepsilon}}_{\text{i}}$
and the inelastic flow $(\dot{\boldsymbol{\varepsilon}}_{\text{ki}},
\dot{\boldsymbol{\varepsilon}}_{\text{di}}, \dot{s}_\text{d})$
which takes place on the microstructural level.

The specific free energy per unit mass is given by
\begin{equation}\label{AddEner}
\psi = \psi_{\text{el}} (\boldsymbol{\varepsilon}_{\text{e}}) +
\psi_{\text{kin}} (\boldsymbol{\varepsilon}_{\text{ke}}) +
\psi_{\text{dis}} (\boldsymbol{\varepsilon}_{\text{de}}) +  \psi_{\text{iso}} (s_\text{e}),
\end{equation}
\begin{equation}\label{concrete1}
\rho \psi_{\text{el}} (\boldsymbol{\varepsilon}_{\text{e}}) =
\frac{k}{2} (\text{tr} \ \boldsymbol{\varepsilon}_{\text{e}})^2 +
\mu \ \| \boldsymbol{\varepsilon}^{\text{D}}_{\text{e}} \|^2, \quad
\rho  \psi_{\text{kin}} (\boldsymbol{\varepsilon}_{\text{ke}}) = \frac{c_{\text{k}}}{2}
 \ \|\boldsymbol{\varepsilon}^{\text{D}}_{\text{ke}} \|^2,
\end{equation}
\begin{equation}\label{concrete2}
\rho \psi_{\text{dis}} (\boldsymbol{\varepsilon}_{\text{de}}) =  \frac{c_{\text{d}}}{2} \
\|\boldsymbol{\varepsilon}^{\text{D}}_{\text{de}} \|^2, \quad
\rho \psi_{\text{iso}}(s_{\text{e}}) = \frac{\gamma}{2} (s_{\text{e}})^2.
\end{equation}
Here, $k, \mu, c_{\text{k}}, c_{\text{d}}, \gamma \geq 0$ are
material parameters; $\rho >0$ stands for the mass density. The
quantity $\psi_{\text{el}} (\boldsymbol{\varepsilon}_{\text{e}})$
stands for the energy stored due to macroscopic elastic
deformations. The remaining part $\psi_{\text{kin}}  +
\psi_{\text{dis}}  +  \psi_{\text{iso}}$ is used to capture the
energy associated with the defects of the crystal
structure.\footnote{Note that $\psi_{\text{kin}} + \psi_{\text{dis}}
+  \psi_{\text{iso}}$ does not necessarily reflect the entire
``defect energy". It is natural to assume that a part of the  ``defect
energy" is not connected to any hardening mechanism
\citep{ShutovIhle}. Alternatively, \cite{Henann} prefer to neglect
the free energy storage $\psi_{\text{iso}}$ which is associated with
the isotropic hardening. Interestingly, \cite{FeigenbaumD2} suggest
that the defect energy is released while the distortion of the yield
surface takes place. Thus, again, somewhat smaller energy storage
will be predicted than by assumption \eqref{AddEner}. The choice
between many alternatives is important for the prediction of the
inelastic dissipation and should be based on relevant experimental
observations.} Next, we postulate the following relations for
stresses, backstresses, and isotropic hardening
\begin{equation}\label{potent}
\boldsymbol{\sigma}= \rho \frac{\partial \psi_{\text{el}}
(\boldsymbol{\varepsilon}_{\text{e}})}{\partial \boldsymbol{\varepsilon}_{\text{e}}}, \
\boldsymbol{X}_{\text{k}} = \rho \frac{\partial \psi_{\text{kin}}
(\boldsymbol{\varepsilon}_{\text{ke}})}{\partial \boldsymbol{\varepsilon}_{\text{ke}}}, \
\boldsymbol{X}_{\text{d}} = \rho \frac{\partial \psi_{\text{dis}}
(\boldsymbol{\varepsilon}_{\text{de}})}{\partial \boldsymbol{\varepsilon}_{\text{de}}}, \
R = \rho \frac{\partial \psi_{\text{iso}}
(s_{\text{e}})}{\partial s_{\text{e}}}.
\end{equation}
Substituting \eqref{concrete1} and \eqref{concrete2} into \eqref{potent} we get
\begin{equation}\label{potent2}
\boldsymbol{\sigma}= k \ \text{tr} (\boldsymbol{\varepsilon}_{\text{e}}) \mathbf{1} +
2 \mu \boldsymbol{\varepsilon}^{\text{D}}_{\text{e}}, \
\boldsymbol{X}_{\text{k}} = c_{\text{k}} \boldsymbol{\varepsilon}^{\text{D}}_{\text{ke}}, \
\boldsymbol{X}_{\text{d}} = c_{\text{d}} \boldsymbol{\varepsilon}^{\text{D}}_{\text{de}}, \
R = \gamma s_{\text{e}}.
\end{equation}

On the one hand, these relations can be motivated by the rheological model from Section 2.
On the other hand, as it will be shown in the following,
relations \eqref{potent} will be sufficient for the thermodynamic consistency of the material model.
It follows immediately from $\eqref{potent2}$ that
$\text{tr}\boldsymbol{X}_{\text{k}} = \text{tr}\boldsymbol{X}_{\text{d}} = 0$.

\begin{figure}\centering
\psfrag{U}[m][][1][0]{$\vec{e}_1$}
\psfrag{EL}[m][][1][0]{$\text{El}(\bar{K}( \cdot , \alpha))$}
\psfrag{T}[m][][1][0]{$\theta$}
\psfrag{FB}[m][][1][0]{$\bar{f}(\vec{y})$}
\psfrag{x}[m][][1][0]{$\vec{y}$}
\psfrag{R}[m][][1][0]{$\mathbb{R}^2$}
\psfrag{NO}[m][][1][0]{$\vec{0}$}
\psfrag{A}[m][][1][0]{$1$}
\psfrag{B}[m][][1][0]{$\sqrt{2/3} (K_0+R)$}
\psfrag{S}[m][][1][0]{Deviatoric stress space}
\psfrag{C}[m][][1][0]{$\boldsymbol{X}_{\text{d}}$}
\psfrag{G}[m][][1][0]{$\boldsymbol{X}_{\text{k}} + \boldsymbol{X}_{\text{d}}$}
\psfrag{W}[m][][1][0]{$\boldsymbol{\sigma}^{\text{D}}_{\text{eff}}$}
\psfrag{Y}[m][][1][0]{$\boldsymbol{\sigma}^{\text{D}}$}
\psfrag{F}[m][][1][0]{$f(\boldsymbol{\sigma})$}
\psfrag{K}[m][][1][0]{$\mathfrak{S}$}
\psfrag{L}[m][][1][0]{$\boldsymbol{\sigma}^{\ast} = \boldsymbol{X}_{\text{k}} + \boldsymbol{X}_{\text{d}}$}
\psfrag{M}[m][][1][0]{$\mathbf{n}$}
\psfrag{J}[m][][1][0]{a)}
\psfrag{H}[m][][1][0]{b)}
\scalebox{0.8}{\includegraphics{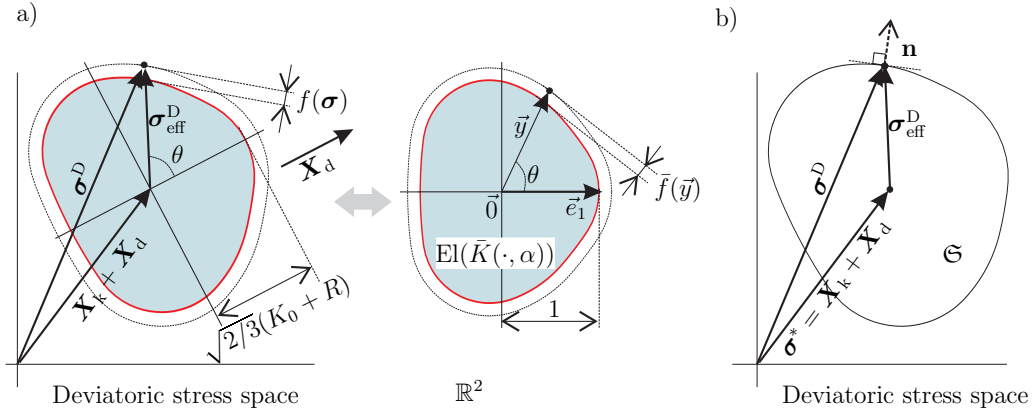}}
\caption{a) Sketch of the yield surface in the deviatoric stress space and definition of the overstress
$f(\boldsymbol{\sigma})$
(cf. \eqref{DefThet}, \eqref{overDe3}). The elastic domain in the stress space is associated
with $\text{El}(\bar{K}( \cdot , \alpha)) \subset \mathbb{R}^2$;
b) Sketch of the proof of the inequality \eqref{Positiv}. \label{fig6}}
\end{figure}

Suppose that the degree of distortion $\alpha$
depends solely on $\| \boldsymbol{X}_{\text{d}} \|$.
A concrete dependence will be specified in the following (cf. \eqref{DistortReal}).
The effective stress tensor and the angle $\theta$ are defined now through
\begin{equation}\label{DefThet}
\boldsymbol{\sigma}_{\text{eff}}  :=   \boldsymbol{\sigma} - \boldsymbol{X}_{\text{k}} - \boldsymbol{X}_{\text{d}}, \quad
\theta := \text{arccos}\Big( \frac{\boldsymbol{\sigma}^{\text{D}}_{\text{eff}} \ : \ \boldsymbol{X}_{\text{d}}}
{\|\boldsymbol{\sigma}^{\text{D}}_{\text{eff}}\| \ \|\boldsymbol{X}_{\text{d}}\|}\Big).
\end{equation}
Note that for $\|\boldsymbol{X}_{\text{d}} \| = 0$ the angle $\theta$ is arbitrary. To be definite, we put
$\theta = 0$ in that case. Further, analogously to \eqref{overDe} and \eqref{overDe2}, we define $\vec{y} \in \mathbb{R}^2$ and the corresponding
overstress $f$ (see Fig. \ref{fig6}a)
\begin{equation}\label{overDe3}
f(\boldsymbol{\sigma}, \boldsymbol{X}_{\text{k}}, \boldsymbol{X}_{\text{d}}, R) =
\tilde{f}(\| \boldsymbol{\sigma}^{\text{D}}_{\text{eff}} \|, \theta, \alpha, R) := \sqrt{\frac{2}{3}} (K_0 + R) \ \bar{f}(\vec{y},\alpha),
\end{equation}
\begin{equation}\label{overDe4}
\vec{y} :=\frac{ \| \boldsymbol{\sigma}^{\text{D}}_{\text{eff}}  \|}{\sqrt{2/3} (K_0+R)}
(\cos(\theta), \sin(\theta)).
\end{equation}
Here, $K_0 >0$ is a fixed material parameter (initial yield stress), and
the function $\bar{f}(\vec{y},\alpha)$ is defined through \eqref{overdef}.

The elastic domain corresponds to stress states with zero overstress $f$.
For a given stress tensor $\boldsymbol{\sigma}$, a non-dimensional vector
$\vec{y} \in \mathbb{R}^2$ must be evaluated according to \eqref{overDe4}.
Observe that the angle between $\vec{y}$ and $\vec{e}_1$ coincides with the angle between
$\boldsymbol{\sigma}^{\text{D}}_{\text{eff}}$ and $\boldsymbol{X}_{\text{d}}$ (see Fig. \ref{fig6}a).
According to \eqref{overDe3}, the stress state $\boldsymbol{\sigma}$ lies within
the elastic domain if and only if $\vec{y} \in \text{El}(\bar{K}( \cdot , \alpha))$.
The origin of the elastic domain corresponds to $\{ \vec{y} = \vec{0} \}$ or, equivalently,
$\{ \boldsymbol{\sigma}^{\text{D}}_{\text{eff}} = \boldsymbol{0} \} =
\{  \boldsymbol{\sigma}^{\text{D}} = \boldsymbol{X}_{\text{k}} + \boldsymbol{X}_{\text{d}} \}$. Next, observe that
the direction of the elastic domain coincides with the direction of $\boldsymbol{X}_{\text{d}}$,
and the size of the elastic domain in that direction equals $\sqrt{\frac{2}{3}} (K_0 + R)$ (see Fig. \ref{fig6}a).

\textbf{Remark 7.} The relations \eqref{overDe3} and \eqref{overDe4} imply that the
resistance to plastic deformation depends on the actual direction of
the loading relative to the recent loading path. This dependence on
the loading direction can be associated with the
activation/deactivation of crystallographic slip planes as well as
mobilization/demobilization of oriented dislocation structures. It
is well known that the gliding of dislocations is obstructed by the
cell walls under monotonic loading ($\theta =0$), but after a strain
path change ($\theta \ne 0$), the loading may drive the dislocations
toward the cell interior \citep{Viatkina2}. In the monograph by
\cite{Viatkina1}, the mechanism of  ``directional remobilisation"
under the strain path change is explained as a remobilisation of
dislocation locks and dipoles which were formed during the previous
loading, in contrast to the statistical remobilisation which is
independent of the strain path change.

The normality flow rule in combination with the Perzyna rule is used
\begin{equation}\label{normFlow}
\dot{\boldsymbol{\varepsilon}}_{\text{i}} =
\lambda_{\text{i}} \frac{ \partial f(\boldsymbol{\sigma}, \boldsymbol{X}_{\text{k}}, \boldsymbol{X}_{\text{d}}, R) }
{\partial \boldsymbol{\sigma}}
\ \text{for} \ f>0, \
\dot{\boldsymbol{\varepsilon}}_{\text{i}} = \boldsymbol 0 \ \text{for} \ f=0;
\quad \lambda_{\text{i}} = \frac{\displaystyle 1}{\displaystyle
\eta}\Big( \frac{\displaystyle 1}{\displaystyle k_0} f \Big)^{m}.
\end{equation}
Here, $ \lambda_{\text{i}} \geq 0$ is an inelastic multiplier which controls the rate
of the inelastic flow. Indeed, since $\| \frac{ \partial \bar{f}(\vec{y},\alpha)}{\partial \vec{y}} \| = 1$
for $\bar{f} >0$, we get (cf. \eqref{App9})
\begin{equation}\label{Odqi}
\Big\| \frac{ \partial f(\boldsymbol{\sigma}, \boldsymbol{X}_{\text{k}}, \boldsymbol{X}_{\text{d}}, R) }
{\partial \boldsymbol{\sigma}} \Big\| =1 \ \text{for} \ f>0, \quad
\| \dot{\boldsymbol{\varepsilon}}_{\text{i}} \| = \lambda_{\text{i}}, \quad
\dot{p} = \lambda_{\text{i}}.
\end{equation}
We emphasize that $k_0$ is not a material parameter, and
we put $k_0 = 1$ MPa.
A concrete algorithm for the evaluation of the derivative
$\frac{ \partial f(\boldsymbol{\sigma})}{\partial \boldsymbol{\sigma}}$ is presented in Appendix B.
Note that the normality rule \eqref{normFlow} implies an incompressible
flow: $\text{tr} \dot{\boldsymbol{\varepsilon}}_{\text{i}} =0$.
In order to take the saturation of the kinematic and distortional hardening
into account, we postulate for the inelastic flows on the microstructural level (cf. \eqref{satur})
\begin{equation}\label{satur2}
\dot{\boldsymbol{\varepsilon}}_{\text{ki}} = \lambda_{\text{i}} \ \varkappa_\text{k} \ \boldsymbol{X}_{\text{k}}, \quad
\dot{\boldsymbol{\varepsilon}}_{\text{di}} = \lambda_{\text{i}} \ \varkappa_\text{d} \ \boldsymbol{X}_{\text{d}}.
\end{equation}
Here, $\varkappa_\text{k}, \varkappa_\text{d} \geq 0$ are material parameters. Recall that
$\text{tr}\boldsymbol{X}_{\text{k}} = \text{tr}\boldsymbol{X}_{\text{d}} = 0$. Thus, the inelastic
flow on the microstructural level is incompressible as well:
$\text{tr} \dot{\boldsymbol{\varepsilon}}_{\text{ki}} = \text{tr} \dot{\boldsymbol{\varepsilon}}_{\text{di}} =0$.
It can be easily shown that for $\| \boldsymbol{X}_{\text{d}} \| |_{t=0} \leq  1/\varkappa_\text{d}$
we have $\| \boldsymbol{X}_{\text{d}} \| \leq 1/\varkappa_\text{d}$.
Thus, analogously to \eqref{DistortRheolo}, we define the distortion parameter $\alpha \in [0,1]$ by
\begin{equation}\label{DistortReal}
\alpha := \varkappa_{\text{d}} \ \|\boldsymbol{X}_{\text{d}}\|.
\end{equation}
For a given deviatoric stress $\boldsymbol{\sigma}^{\text{D}}$ consider a convex set
$\mathfrak{S} = \{\boldsymbol{\sigma}^{\ast} : \ \text{tr} \boldsymbol{\sigma}^{\ast} =0,
\quad  f(\boldsymbol{\sigma}^{\ast}) \leq f(\boldsymbol{\sigma}^{\text{D}})\}$.
The gradient $\frac{ \partial f(\boldsymbol{\sigma}, \boldsymbol{X}_{\text{k}}, \boldsymbol{X}_{\text{d}}, R) }
{\partial \boldsymbol{\sigma}}$ coincides with the unit outward normal $\mathbf{n}$ to the boundary of $\mathfrak{S}$
at $\boldsymbol{\sigma}^{\text{D}}$.
Moreover, the state $\boldsymbol{\sigma}^{\ast} = \boldsymbol{X}_{\text{k}} + \boldsymbol{X}_{\text{d}}$
lies within $\mathfrak{S}$
(cf. Fig. \ref{fig6}b). Due to the convexity of $\mathfrak{S}$ we have
\begin{equation}\label{Positiv}
\boldsymbol{\sigma}^{\text{D}}_{\text{eff}} : \dot{\boldsymbol{\varepsilon}}_{\text{i}}
\stackrel{\eqref{normFlow}}{=}
\lambda_{\text{i}} (\boldsymbol{\sigma}^{\text{D}} -  \boldsymbol{\sigma}^{\ast}) : \mathbf{n} \geq 0.
\end{equation}
Having this inequality in mind we formulate the evolution equations for the internal variable
$s$ and its dissipative part $s_{\text{d}}$:
\begin{equation}\label{Odqvist}
\dot{s} = \frac{\boldsymbol{\sigma}_{\text{eff}} : \dot{\boldsymbol{\varepsilon}}_{\text{i}}}{K_0 + R}
\stackrel{\text{tr} \dot{\boldsymbol{\varepsilon}}_{\text{i}} =0}{=}
\frac{\boldsymbol{\sigma}^{\text{D}}_{\text{eff}} : \dot{\boldsymbol{\varepsilon}}_{\text{i}}}{K_0 + R}, \quad
\quad
\dot{s}_{\text{d}} = \frac{\beta}{\gamma} \dot{s} R,
\end{equation}
where $\beta \geq 0$ is a material parameter controlling the saturation of the isotropic hardening.
It follows from \eqref{Positiv} that $\dot{s} \geq 0$. Thus, similar to the inelastic
arc-length, the variable $s$ increases monotonically.
Note that in the case of proportional monotonic loading we have $\theta \approx 0$. Thus,
$\boldsymbol{\sigma}^{\text{D}}_{\text{eff}} : \dot{\boldsymbol{\varepsilon}}_{\text{i}} \approx
\| \boldsymbol{\sigma}^{\text{D}}_{\text{eff}} \| \ \| \dot{\boldsymbol{\varepsilon}}_{\text{i}} \| =
\| \boldsymbol{\sigma}^{\text{D}}_{\text{eff}} \| \ \lambda_{\text{i}}$.
Moreover, for slow loading we get $f \ll (K_0+R)$. Therefore,
$\| \boldsymbol{\sigma}^{\text{D}}_{\text{eff}} \| \approx \sqrt{2/3} \ (K_0 + R)$.
Thus, under quasistatic proportional loading, the parameter $s$ evolves similar to the
inelastic arc-length: $\dot{s} \approx  \sqrt{2/3} \
\| \dot{\boldsymbol{\varepsilon}}_{\text{i}} \| = \sqrt{2/3} \lambda_{\text{i}}$.
Under general loading conditions, although, the evolution of $s$ depends not only
on the rate of the plastic flow, but also on its direction.

Finally, the system of constitutive equations is closed
by initial conditions imposed on the strain-like internal variables
\begin{equation*}\label{InitCond}
\boldsymbol{\varepsilon}_{\text{i}} |_{t=0}  = \boldsymbol{\varepsilon}^0_{\text{i}}, \
\boldsymbol{\varepsilon}_{\text{ki}} |_{t=0}  = \boldsymbol{\varepsilon}^0_{\text{ki}}, \
\boldsymbol{\varepsilon}_{\text{di}} |_{t=0}  = \boldsymbol{\varepsilon}^0_{\text{di}}, \
s|_{t=0} = s^0, \ s_{\text{d}} |_{t=0} = s^0_{\text{d}}.
\end{equation*}
We suppose
$\text{tr} \boldsymbol{\varepsilon}^0_{\text{i}} =
\text{tr} \boldsymbol{\varepsilon}^0_{\text{ki}} =
\text{tr} \boldsymbol{\varepsilon}^0_{\text{di}} =0$.
If the undeformed state is assumed to be stress free at $t=0$, then
$\boldsymbol{\varepsilon}^0_{\text{i}} = \boldsymbol 0$.
The quantities $\boldsymbol{\varepsilon}^0_{\text{ki}}$ and
$\boldsymbol{\varepsilon}^0_{\text{di}}$ can be used to capture the initial plastic
 anisotropy of the material.\footnote{This is equivalent to
 the introduction of initial backstresses.} In particular,
the yield condition at $t=0$ does not have to coincide with the Huber-Mises criterium.

%\textbf{Remark}. For large $\varkappa_d$ and $\varkappa_k$,
%the yield surface can be
%distorted without being translated away from the zero point in the
%stress space.

\subsection{Proof of thermodynamic consistency}

Let us consider the Clausius-Duhem inequality in the form (see, for example, \cite{Haupt})
\begin{equation}\label{CDU}
\delta_{\text{i}} := \frac{1}{\rho} \boldsymbol{\sigma} : \dot{\boldsymbol{\varepsilon}} - \dot{\psi} \geq 0.
\end{equation}
Taking the kinematic relations \eqref{AddDec} into account, we rewrite the stress power as follows
\begin{multline}\label{CDU2}
\boldsymbol{\sigma} : \dot{\boldsymbol{\varepsilon}} = \\
 \boldsymbol{\sigma} : ( \dot{\boldsymbol{\varepsilon}_{\text{e}}} +
\dot{\boldsymbol{\varepsilon}_{\text{i}}})
- \boldsymbol{X}_{\text{k}} : \dot{\boldsymbol{\varepsilon}_{\text{i}}} +
\boldsymbol{X}_{\text{k}} : (\dot{\boldsymbol{\varepsilon}_{\text{ki}}} +
\dot{\boldsymbol{\varepsilon}_{\text{ke}}}) -
\boldsymbol{X}_{\text{d}} : \dot{\boldsymbol{\varepsilon}_{\text{i}}} +
\boldsymbol{X}_{\text{d}} : ( \dot{\boldsymbol{\varepsilon}_{\text{di}}} +
\dot{\boldsymbol{\varepsilon}_{\text{de}}}).
\end{multline}
Moreover, differentiating \eqref{AddEner}, we get for the time derivative of the free energy
\begin{equation}\label{CDU3}
\dot{\psi} = \frac{\partial \psi_{\text{el}}
(\boldsymbol{\varepsilon}_{\text{e}})}{\partial \boldsymbol{\varepsilon}_{\text{e}}} :
\dot{\boldsymbol{\varepsilon}}_{\text{e}} +
\frac{\partial \psi_{\text{kin}}
(\boldsymbol{\varepsilon}_{\text{ke}})}{\partial \boldsymbol{\varepsilon}_{\text{ke}}} :
\dot{\boldsymbol{\varepsilon}}_{\text{ke}} +
\frac{\partial \psi_{\text{dis}}
(\boldsymbol{\varepsilon}_{\text{de}})}{\partial \boldsymbol{\varepsilon}_{\text{de}}} :
\dot{\boldsymbol{\varepsilon}}_{\text{de}} +
\frac{\partial \psi_{\text{iso}}
(s_{\text{e}})}{\partial s_{\text{e}}} \dot{s}_{\text{e}}.
\end{equation}
Substituting \eqref{CDU2} and \eqref{CDU3} into \eqref{CDU} and taking
the potential relations \eqref{potent} into account, we obtain
the Clausius-Duhem inequality in the following form
\begin{equation*}\label{CDU4}
\rho \ \delta_{\text{i}} = \big( \boldsymbol{\sigma}_{\text{eff}} : \dot{\boldsymbol{\varepsilon}}_{\text{i}} -
R \ \dot{s} \big) + \boldsymbol{X}_{\text{k}} : \dot{\boldsymbol{\varepsilon}}_{\text{ki}} +
\boldsymbol{X}_{\text{d}} : \dot{\boldsymbol{\varepsilon}}_{\text{di}} + R \dot{s}_{\text{d}} \geq 0.
\end{equation*}
It follows immediately from \eqref{satur2} and $\eqref{Odqvist}_2$ that
$\boldsymbol{X}_{\text{k}} : \dot{\boldsymbol{\varepsilon}}_{\text{ki}} \geq 0$,
$\boldsymbol{X}_{\text{d}} : \dot{\boldsymbol{\varepsilon}}_{\text{di}} \geq 0$, and
$R \dot{s}_{\text{d}} \geq 0$. In order to prove the thermodynamic consistency of the material
model it remains to show that $\boldsymbol{\sigma}_{\text{eff}} : \dot{\boldsymbol{\varepsilon}}_{\text{i}} -
R \ \dot{s} \geq 0$. Indeed,
\begin{equation*}\label{CDU5}
\boldsymbol{\sigma}_{\text{eff}} : \dot{\boldsymbol{\varepsilon}}_{\text{i}} -
R \ \dot{s}  \stackrel{\eqref{Odqvist}_1}{=} \boldsymbol{\sigma}^{\text{D}}_{\text{eff}} : \dot{\boldsymbol{\varepsilon}}_{\text{i}} \
(1 - R/(K_0+R)) = \boldsymbol{\sigma}^{\text{D}}_{\text{eff}} : \dot{\boldsymbol{\varepsilon}}_{\text{i}} \
(K_0/(K_0+R)) \stackrel{\eqref{Positiv}}{\geq} 0.
\end{equation*}
The thermodynamic consistency of the material model is thus proved.

\textbf{Remark 8.}
Note that the proof of the thermodynamic consistency
is essentially based on the inequality
$\boldsymbol{\sigma}^{\text{D}}_{\text{eff}} : \dot{\boldsymbol{\varepsilon}}_{\text{i}} \geq 0$.
\emph{Any flow rule} which governs $\dot{\boldsymbol{\varepsilon}}_{\text{i}}$ and
complies with this inequality
would yield a thermodynamically consistent material model, as well.
For instance, the radial flow rule can be considered as a simplified alternative
to the normality rule $\eqref{normFlow}_1$
\begin{equation}\label{normFlow342}
\dot{\boldsymbol{\varepsilon}}_{\text{i}} = \lambda_{\text{i}} \boldsymbol{R}_{\text{eff}}, \quad
\boldsymbol{R}_{\text{eff}}:=\frac{ \boldsymbol{\sigma}^{\text{D}}_{\text{eff}} }
{\| \boldsymbol{\sigma}^{\text{D}}_{\text{eff}} \|  }.
\end{equation}

\subsection{Identification of material parameters}

The material model contains 11 material parameters and a
material function $\bar{K}^{\text{sat}}(\theta)$. Let us discuss
the identification of these quantities. First,
the elasticity parameters $k$ and $\mu$ can be determined
basing on the experimental data for elastic deformations. Next,
the initial yield stress $K_0$ can be calibrated using the graphical
method from a quasistatic uniaxial tension test. The viscosity
parameters $\eta$ and $m$ of the Perzyna law are typically identified
using a series of tests under monotonic loading with
different loading rates.
Further, the material function
$\bar{K}^{\text{sat}}(\theta)$ is uniquely determined for the
fixed $\bar{K}^{\text{sat}}(\pi)$ if
the form of the saturated yield surface is known
(for details see Remark 6).
For simplicity, one may assume $\bar{K}^{\text{sat}}(\pi)=1$.
In that case, the parameters of isotropic hardening ($\gamma$ and $\beta$)
can be identified using the information about how the size of the elastic domain
evolves under monotonic loading. Finally, it remains to identify two parameters of
kinematic hardening ($c_{\text{k}}$ and $\varkappa_{\text{k}}$)
and two distortional parameters ($c_{\text{d}}$ and $\varkappa_{\text{d}}$).
This can be done by minimization of a
least-squares functional which represents the discrepancy
between measurements data and corresponding
model predictions. Experimental measurements related to non-proportional loading
are necessary in order to obtain a reliable identification procedure.
Some regularization techniques
can be used to reduce the correlation among the parameters and to
reduce the probability of getting trapped in local minima \citep{Shutov3}.

The success of the identification procedure
depends on the quality of initial approximation chosen for
the unknown parameters $c_{\text{k}}$,
$c_{\text{d}}$, $\varkappa_{\text{k}}$, and $\varkappa_{\text{d}}$.
The order of magnitude of these parameters can be estimated basing on
the following considerations: The upper bounds for $\|\boldsymbol{X}_{\text{k}}\|$ and $\|\boldsymbol{X}_{\text{d}}\|$
are given by $\varkappa_{\text{k}}^{-1}$ and $\varkappa_{\text{d}}^{-1}$, respectively.
At the same time, the increment of the inelastic arc-length which corresponds to
the saturation of kinematic and distortional hardening under proportional loading
is proportional to $(c_{\text{k}} \varkappa_{\text{k}})^{-1}$
and $(c_{\text{d}} \varkappa_{\text{d}})^{-1}$, respectively.

%The germ center lies exactly at $Xk+Xd$. Therefore, it translates kinematically in the stress space.

\section{Numerical computations}

In this study, for simplicity, the evolution equations
\eqref{normFlow}, \eqref{satur2}, and \eqref{Odqvist} are
integrated numerically using explicit time-stepping scheme.
If rate-independent material response is to be simulated,
a viscous regularization with fictitious small
viscosity $\eta > 0$ can be used.\footnote{Note that
such viscous regularization allows to smoothen the
sharp transition between elastic and plastic regions.}

In this section we validate the predictive capabilities of the material model.
Toward that end, we consider experimental data of \cite{KhanII} obtained for
a very high work hardening aluminum alloy - annealed 1100 Al.
The yield points were identified experimentally under combined tension-torsion
of thin-walled tubular specimens using a small proof strain.
In order to simulate the deformation of thin-walled tubular specimen
we compute the stress response at a single material point.
Consider a Cartesian coordinate system
such that its basis vectors $\mathbf{e}_1$, $\mathbf{e}_2$, and $\mathbf{e}_3$
are oriented along the local axial, hoop, and radial directions, respectively.
The stress state can be idealized approximately as a
special case of the plane stress:
\begin{equation}\label{planeStress}
\boldsymbol{\sigma} = \sigma_{11} \mathbf{e}_1 \otimes \mathbf{e}_1 +
\sigma_{12} (\mathbf{e}_1 \otimes \mathbf{e}_2 + \mathbf{e}_2 \otimes \mathbf{e}_1).
\end{equation}

\begin{figure}\centering
\psfrag{R}[m][][1][0]{$\sqrt{3} \ \sigma_{12}$ [MPa]}
\psfrag{S}[m][][1][0]{$\sigma_{11}$ [MPa]}
\psfrag{J}[m][][1][0]{a)}
\psfrag{H}[m][][1][0]{b)}
\scalebox{0.8}{\includegraphics{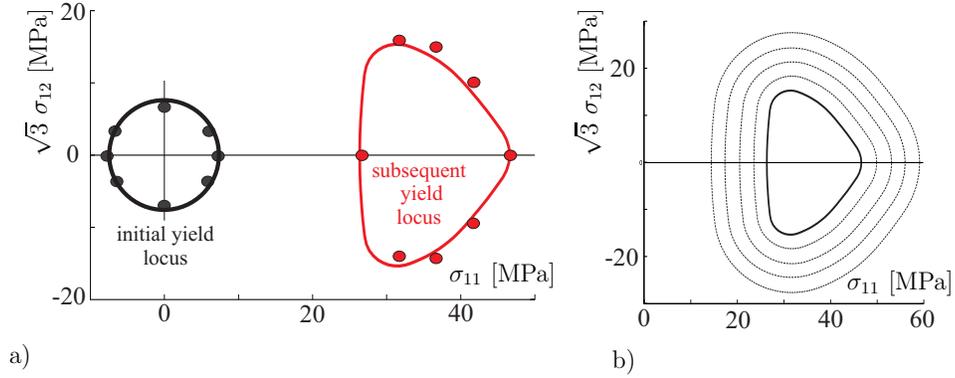}}
\caption{a) Experimental data for annealed 1100 aluminum alloy
\citep{KhanII} and corresponding simulation results,
b) Lines of constant overstress $f$ corresponding to
plastic anisotropy induced by 2\% prestrain in the axial direction. \label{fig7}}
\end{figure}

Here, $\sigma_{11}$ and $\sigma_{12}$ are associated to
the axial and torsional loading, respectively.
The measurement results are
represented in the $(\sigma_{11} , \sqrt{3} \sigma_{12})$-space in Fig. \ref{fig7}a.
The initial yield surface
can be approximated with sufficient accuracy using the
conventional Huber-Mises yield condition.
Therefore, an initial plastic
isotropy will be assumed during the material modeling.
The yield points which were determined after 2\% axial prestrain
are depicted in Fig. \ref{fig7}a, as well. The material parameters used
to simulate the material response are summarized in Table \ref{tab1}.\footnote{It
is not the aim of the current study to identify
the material parameters corresponding to the 1100 aluminum alloy.
Instead, we validate the material model
by the qualitative description of the real experimental data.}
Moreover, due to the initial isotropy we consider the initial conditions as follows:
\begin{equation*}\label{incon}
\boldsymbol{\varepsilon}^0_{\text{i}} = \boldsymbol{\varepsilon}^0_{\text{di}}
= \boldsymbol{\varepsilon}^0_{\text{ki}} = \mathbf{0}, \quad s^0=s^0_{\text{d}}=0.
\end{equation*}
The smooth function
$\bar{K}^{\text{sat}}(\theta)$ which is needed to compute
the overstress $\bar{f}(\vec{y})$ corresponds to
the saturated form shown in Fig. \ref{fig3}c.
Due to the axial prestrain, the distortion parameter
$\alpha$ ranges from $0$ up to $0.999993$, which corresponds
to the (almost) saturated distortional hardening.
As it is shown in Fig. \ref{fig7}a, the yield locus
undergoes the isotropic expansion, kinematic translation and distortion.
The surfaces of constant overstress are depicted in Fig. \ref{fig7}b.
In accordance with the modeling assumptions, these surfaces are slightly less
distorted than the corresponding yield surface.

\textbf{Remark 9.} Note that the form of the yield surface
in the $(\sigma_{11} , \sqrt{3} \sigma_{12})$-space
coincides with the boundary of $\text{El}(\bar{K}( \cdot , \alpha))$.
This is due to the well-known fact that the scalar product of two symmetric
tensors $\boldsymbol{\sigma}^{\text{D}}_{I}$ and $\boldsymbol{\sigma}^{\text{D}}_{II}$
(where $\boldsymbol{\sigma}_{I}$ and $\boldsymbol{\sigma}_{II}$ comply with \eqref{planeStress}) corresponds
to the product of two vectors $\vec{\sigma}_{I}, \vec{\sigma}_{II} \in \mathbb{R}^2$, defined by
$\vec{\sigma}_{I} := (\sigma^{I}_{11}, \sqrt{3} \sigma^{I}_{12})$,
$\vec{\sigma}_{II} := (\sigma^{II}_{11}, \sqrt{3} \sigma^{II}_{12})$. More precisely
\begin{equation*}\label{isomorph}
\boldsymbol{\sigma}^{\text{D}}_{I} : \boldsymbol{\sigma}^{\text{D}}_{II} =
\frac{2}{3} (\sigma^{I}_{11} \sigma^{II}_{11} + 3 \sigma^{I}_{12} \sigma^{II}_{12})
= \frac{2}{3} \vec{\sigma}_{I} \cdot \vec{\sigma}_{II}.
\end{equation*}
Therefore, in the context of \eqref{planeStress},
the angle between two deviatoric stress-states coincides with the angle
between two corresponding vectors in the $(\sigma_{11}, \sqrt{3} \sigma_{12})$-space:
\begin{equation*}\label{isomorph2}
\arccos\Big(\frac{\boldsymbol{\sigma}^{\text{D}}_{I} : \boldsymbol{\sigma}^{\text{D}}_{II} }
{ \|\boldsymbol{\sigma}^{\text{D}}_{I}\| \ \| \boldsymbol{\sigma}^{\text{D}}_{II} \| }\Big) =
\arccos\Big( \frac{\vec{\sigma}_{I} \cdot \vec{\sigma}_{II} }
{ \| \vec{\sigma}_{I} \| \ \| \vec{\sigma}_{II} \| }\Big).
\end{equation*}
Thus, if the form of the saturated yield surface is determined experimentally in
the $(\sigma_{11} , \sqrt{3} \sigma_{12})$-space, it can be used
to identify the smooth function $\bar{K}^{\text{sat}}(\theta)$.

\begin{table}[h]
\caption{Material parameters used to validate the material model: the viscosity
effects are neglected.}
\begin{tabular}{| l l l l l |}
\hline
$k$ [MPa] & $\mu$ [MPa]  &  $c_{\text{k}}$  [MPa]  & $c_{\text{d}}$  [MPa] &  $\gamma$ [MPa] \\ \hline
69000 & 26000  &  1010 & 5000 &  245   \\ \hline
\end{tabular} \\
\begin{tabular}{|l l l l l l |}
\hline
$K_0$ [MPa] & $m$ [-] & $\eta$ [$\text{s}^{-1}$] & $\varkappa_{\text{k}}$ [$\text{MPa}^{-1}$] & $\varkappa_{\text{d}}$ [$\text{MPa}^{-1}$] & $\beta$ [-]   \\ \hline
7.4       & 1     & 0                                      &  0.02                                      &             0.1      & 35            \\ \hline
\end{tabular}
\label{tab1}
\end{table}

\begin{figure}\centering
\psfrag{R}[m][][1][0]{$\sqrt{3} \ \sigma_{12}$ [MPa]}
\psfrag{S}[m][][1][0]{$\sigma_{22}$ [MPa]}
\psfrag{K}[m][][1][0]{$\sigma_{11}=0$ [MPa]}
\psfrag{L}[m][][1][0]{$\sigma_{11}=5$ [MPa]}
\psfrag{M}[m][][1][0]{$\sigma_{11}=10$ [MPa]}
\psfrag{N}[m][][1][0]{$\sigma_{11}=15$ [MPa]}
\psfrag{O}[m][][1][0]{$\sigma_{11}=0$ [MPa]}
\psfrag{P}[m][][1][0]{$\sigma_{11}=-5$ [MPa]}
\psfrag{Q}[m][][1][0]{$\sigma_{11}=-10$ [MPa]}
\psfrag{X}[m][][1][0]{$\sigma_{11}=-15$ [MPa]}
\psfrag{J}[m][][1][0]{a)}
\psfrag{H}[m][][1][0]{b)}
\scalebox{0.7}{\includegraphics{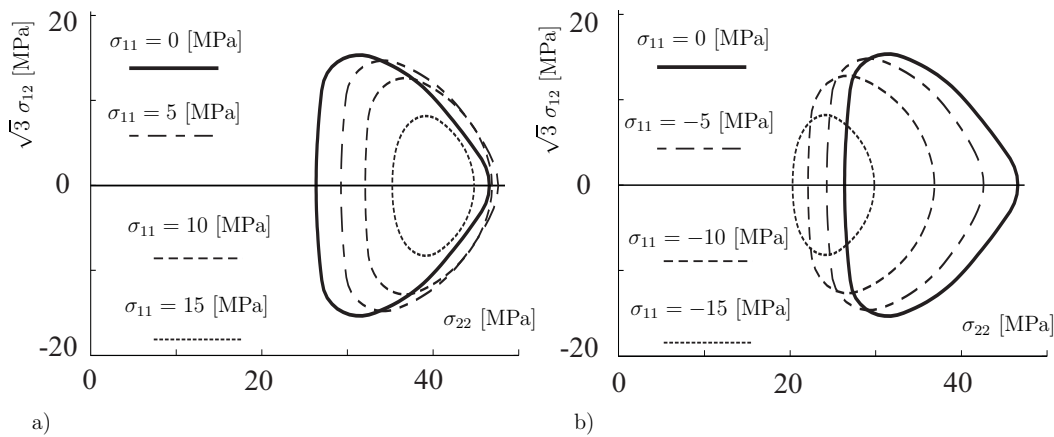}}
\caption{Influence of the axial stress $\sigma_{11}$
on the form of the yield surface after 2\% prestrain in the hoop direction:
a) Positive axial stresses;
b) Negative axial stresses. \label{fig8}}
\end{figure}

If, additionally, the specimen can be loaded by an internal pressure,
the hoop stress $\sigma_{22}$ must be considered, as well:
\begin{equation*}\label{TriaxStress}
\boldsymbol{\sigma} = \sigma_{11} \mathbf{e}_1 \otimes \mathbf{e}_1 + \sigma_{22} \mathbf{e}_2 \otimes \mathbf{e}_2 +
\sigma_{12} (\mathbf{e}_1 \otimes \mathbf{e}_2 + \mathbf{e}_2 \otimes \mathbf{e}_1), \ \sigma_{22} \geq 0.
\end{equation*}
We simulate the evolution of the yield surface in the
process as follows. Starting from the same isotropic initial state,
a 2\% prestrain is prescribed in the hoop direction. Thus, a similar
plastic anisotropy is introduced, as in the previous case.
The form of the yield surface for a fixed axial stress $\sigma_{11}$
is then represented in the
$(\sigma_{22} , \sqrt{3} \sigma_{12})$-space.
As it can be seen in Fig. \ref{fig8}, the form and the size of the yield loci
for $\sigma_{11}=5$ MPa, $\sigma_{11}=10$ MPa, and $\sigma_{11}=15$ MPa
are the same as for $\sigma_{11}=-5$ MPa, $\sigma_{11}=-10$ MPa,
and $\sigma_{11}=-15$ MPa, respectively. Similar to the conventional
Huber-Mises yield condition, the material yields at larger $\sigma_{22}$
stresses for positive $\sigma_{11}$ than for negative $\sigma_{11}$.

%To visualize the evolution of the yield locus let us consider its
%shape in the principal stress space $(\sigma_1 - \sigma_2)$.

%Qualitatively similar evolution of the yield surface was
%experimentally observed by \cite{Rees} for copper-alloys.

\section{Conclusion}

A new material model of metal viscoplasticity with an extremely simple structure
is presented in the current study.
The main modeling assumptions are visualized with the help of
a new two-dimensional rheological model. Only second-rank backstress-like
tensors are used to capture the path dependent
evolution of the plastic anisotropy. Thus, mathematically,
the model is not much more complicated
than the classical model of \cite{Chaboche1, Chaboche2}. At the same time,
the proposed technique possess considerable generality and flexibility.
No specific form of the saturated yield locus is considered, since \emph{any} smooth convex
yield locus can be captured.

An important ingredient of the material modeling is the interpolation between
the initial intact yield surface and the saturated one. The interpolation rule
proposed in the current study \emph{ensures the convexity} of the yield surface at any
stage of hardening. It is shown that this interpolation rule
allows to obtain a thermodynamically consistent material model.

The model contains 6 hardening parameters with 2 parameters per hardening type. These
parameters posses a clear mechanical interpretation and can be identified experimentally.

\section*{Acknowledgement}

This research was supported by German Research Foundation (DFG)
within SFB 692.

\section*{Appendix A}

\begin{figure}\centering
\psfrag{o}[m][][1][0]{$\vec{0}$}
\psfrag{A}[m][][1][0]{$\vec{y}^{\ 0}$}
\psfrag{B}[m][][1][0]{$\vec{y}^{\ 1}$}
\psfrag{C}[m][][1][0]{$\vec{y}_c^{\ 1}$}
\psfrag{D}[m][][1][0]{$\vec{y}_c^{\ i}$}
\psfrag{S}[m][][1][0]{$\alpha \vec{y}_c^{\ i}$}
\psfrag{E}[m][][1][0]{$\vec{y}^{\ i-1}$}
\psfrag{J}[m][][1][0]{$\alpha \vec{y}^{\ i-1}$}
\psfrag{F}[m][][1][0]{$\vec{y}^{\ i}$}
\psfrag{U}[m][][1][0]{$\alpha \vec{y}^{\ i}$}
\psfrag{G}[m][][1][0]{$\vec{y}^{\ N}$}
\psfrag{H}[m][][1][0]{$\vec{n}^{\ i-1}$}
\psfrag{L}[m][][1][0]{$\vec{n}^{\ i}$}
\psfrag{K}[m][][1][0]{$\vec{t}^{\ i-1}$}
\psfrag{M}[m][][1][0]{$\vec{t}^{\ i}$}
\psfrag{MM}[m][][1][0]{$-\vec{t}^{\ i}$}
\psfrag{X}[m][][1][0]{$\vec{y}$}
\psfrag{P}[m][][1][0]{$\mathbf{Q} \cdot (\vec{y}^{\ i} - \vec{y}^{\ i-1})$}
\psfrag{T}[m][][1][0]{$\text{El}^{\text{sat}}$}
\psfrag{Q}[m][][1][0]{$\alpha r^i$}
\psfrag{R}[m][][1][0]{$\alpha \text{El}^{\text{sat}}$}
\psfrag{Y}[m][][1][0]{a)}
\psfrag{Z}[m][][1][0]{b)}
\scalebox{0.8}{\includegraphics{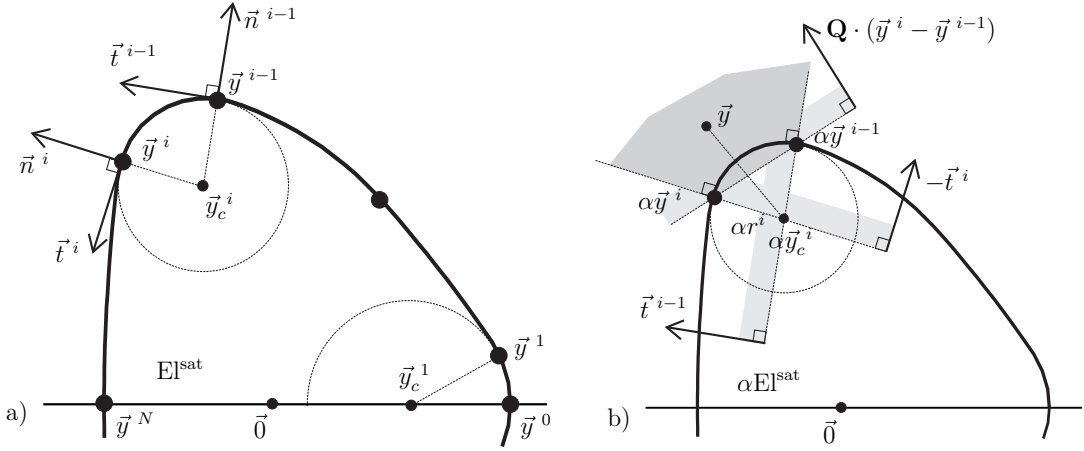}}
\caption{a) A smooth boundary of $\text{El}^{\text{sat}}$ is represented
by a sequence of circular arcs. Each circle is characterized by its center $\vec{y}_c^{\ i}$ and
radius $r^i$;
b) For a given $\vec{y}$,
the inequalities \eqref{ine1}, \eqref{ine2} are
satisfied within the shaded region.
For the suitable $i$, the distance to $\alpha \text{El}^{\text{sat}}$
is computed as
$\mathcal{D} \big(\vec{y}, \alpha \text{El}^{\text{sat}} \big) =
\| \vec{y} - \alpha  \vec{y}_c^{\ i} \| - \alpha r^i$.
\label{fig9}}
\end{figure}

Let us discuss the numerical computation of
$\mathcal{D} \big(\vec{y}, \alpha \text{El}^{\text{sat}} \big)$.
Suppose that the upper half of the boundary of $\text{El}^{\text{sat}}$
is given by $N$ circular arcs connecting $\vec{y}^{\ i-1}$ and $\vec{y}^{\ i}$
for $i \in \{1,...,N\}$ as shown in Fig. \ref{fig9}a. The outward normal
and the tangent at $\vec{y}^{\ i}$ will be denoted
by $\vec{n}^{\ i}$ and $\vec{t}^{\ i}$, respectively. The orientation of
the tangent is chosen in such way that the pair $\{ \vec{n}^{\ i}, \vec{t}^{\ i} \}$
forms a right-handed corner.
In particular, we have
\begin{equation*}
\vec{y}^{\ 0} = (1,0), \quad
\vec{t}^{\ 0} = (0,1), \quad
\vec{t}^{\ N} = (0,-1).
\end{equation*}
For each arc connecting $\vec{y}^{\ i-1}$ and $\vec{y}^{\ i}$ consider
its center $\vec{y}_c^{\ i}$ and its radius $r^{i}$.
In order to make sure that the boundary of $\text{El}^{\text{sat}}$ is smooth, we require
\begin{equation*}
\frac{\vec{y}^{\ 0} - \vec{y}_c^{\ 1}}{ r^{1}} = \vec{n}^{\ 0}, \quad
\frac{\vec{y}^{\ N} - \vec{y}_c^{\ N}}{ r^{N}} = \vec{n}^{\ N};
\end{equation*}
\begin{equation*}
\frac{\vec{y}^{\ i} - \vec{y}_c^{\ i}}{ r^{i}} =
\frac{\vec{y}^{\ i} - \vec{y}_c^{\ i+1}}{ r^{i+1}} = \vec{n}^{\ i} \ \text{for all} \
i \in \{1,...,N-1 \}.
\end{equation*}

Let $\alpha \in [0,1]$ and $\vec{y} = \| \vec{y} \| (\cos(\theta),\sin(\theta))$ be
given such that
$\theta \in [0, \pi]$.
The algorithm used to compute the distance
$\mathcal{D} \big(\vec{y}, \alpha \text{El}^{\text{sat}} \big)$ is as follows.
\begin{itemize}
\item Check the inclusion: If $\| \vec{y} \| \leq \alpha \bar{K}^{\text{sat}}(\theta)$
then $\mathcal{D} \big(\vec{y}, \alpha \text{El}^{\text{sat}} \big)=0$.
\item Otherwise, find the corresponding arc: Find $i \in \{1,...,N\}$ such that (cf. Fig. \ref{fig9}b)
\begin{equation}\label{ine1}
(\vec{y} - \alpha \vec{y}^{\ i-1}) \cdot \vec{t}^{\ i-1} \geq 0, \quad
(\vec{y} - \alpha \vec{y}^{\ i}) \cdot (- \vec{t}^{\ i} ) \geq 0, \
\end{equation}
\begin{equation}\label{ine2}
(\vec{y} - \alpha \vec{y}^{\ i-1}) \cdot  \mathbf{Q} \cdot (\vec{y}^{\ i} - \vec{y}^{\ i-1}) \geq 0, \quad
\mathbf{Q} := \vec{e}_1 \otimes \vec{e}_2 - \vec{e}_2 \otimes \vec{e}_1.
\end{equation}
\item The distance is then given by
\begin{equation*}
\mathcal{D} \big(\vec{y}, \alpha \text{El}^{\text{sat}} \big) = \| \vec{y} - \alpha  \vec{y}_c^{\ i} \| - \alpha r^i.
\end{equation*}
Moreover, the outward unit normal is given by
\begin{equation}\label{ApA8}
\vec{n} = \frac{\partial \bar{f}\big(\vec{y},\alpha \big) }
{\partial \vec{y}} = \frac{\partial  \mathcal{D} \big(\vec{y}, \alpha \text{El}^{\text{sat}} \big)}{\partial
\vec{y}}  = \frac{\vec{y} - \alpha  \vec{y}_c^{\ i}}{\| \vec{y} - \alpha  \vec{y}_c^{\ i} \|}.
\end{equation}
\end{itemize}

\section*{Appendix B}

Let us discuss the computation of the derivative
$\frac{ \partial f(\boldsymbol{\sigma})}{\partial \boldsymbol{\sigma}}$,
which enters the formulation of
the normality rule \eqref{normFlow}. Suppose $f>0$. Thus,
$\boldsymbol{\sigma}^{\text{D}}_{\text{eff}} \neq \boldsymbol0$ and
the radial direction
$\boldsymbol{R}_{\text{eff}} := \boldsymbol{\sigma}^{\text{D}}_{\text{eff}} / \| \boldsymbol{\sigma}^{\text{D}}_{\text{eff}} \|$
is well defined.
Recall that the hardening variables $\boldsymbol{X}_{\text{k}}$, $\boldsymbol{X}_{\text{d}}$, and
$R$ are to be held constant during differentiation.
Having this in mind, we get
\begin{equation}\label{App1}
\frac{ \partial  \| \boldsymbol{\sigma}^{\text{D}}_{\text{eff}} \|}
{\partial \boldsymbol{\sigma}} =\frac{\boldsymbol{\sigma}^{\text{D}}_{\text{eff}}}
{\| \boldsymbol{\sigma}^{\text{D}}_{\text{eff}} \|} = \boldsymbol{R}_{\text{eff}}.
\end{equation}
Next, taking into account that $\frac{d \ \text{arccos}(\Phi)}{d \Phi} = - ( \sin (\text{arccos}(\Phi)) )^{-1}$,
we obtain for $\theta \neq 0$
\begin{equation}\label{App2}
\frac{ \partial  \theta }
{\partial \boldsymbol{\sigma}} =
\frac{ \partial \  \text{arccos}\Big( \frac{\boldsymbol{\sigma}^{\text{D}}_{\text{eff}} \ : \ \boldsymbol{X}_{\text{d}}}
{\|\boldsymbol{\sigma}^{\text{D}}_{\text{eff}}\| \ \|\boldsymbol{X}_{\text{d}}\|}\Big) }
{\partial \boldsymbol{\sigma}} =
-\frac{1}{\sin \theta} \frac{ \partial  \Big( \frac{\boldsymbol{\sigma}^{\text{D}}_{\text{eff}} \ : \ \boldsymbol{X}_{\text{d}}}
{\|\boldsymbol{\sigma}^{\text{D}}_{\text{eff}}\| \ \|\boldsymbol{X}_{\text{d}}\|}\Big)}
{\partial \boldsymbol{\sigma}}.
\end{equation}
Moreover, since $\theta \in [0, \pi]$ is the angle between $\boldsymbol{X}_{\text{d}}$ and
$\boldsymbol{R}_{\text{eff}}$, we get
\begin{equation*}
\sin \theta = \frac{\| \boldsymbol{X}_{\text{d}} -
(\boldsymbol{X}_{\text{d}}:\boldsymbol{R}_{\text{eff}} ) \boldsymbol{R}_{\text{eff}} \|}
{\|\boldsymbol{X}_{\text{d}}\|}.
\end{equation*}
Substituting this result into \eqref{App2} and taking into account
that $ \frac{\partial ( \boldsymbol{\sigma}^{\text{D}}_{\text{eff}} \ : \ \boldsymbol{X}_{\text{d}}) }
{\partial \boldsymbol{\sigma}} = \boldsymbol{X}_{\text{d}}$ we get for $\theta \neq 0$
\begin{equation}\label{App4}
\frac{ \partial  \theta }
{\partial \boldsymbol{\sigma}} = -
\frac{1}{\|\boldsymbol{\sigma}^{\text{D}}_{\text{eff}}\|} \frac{\boldsymbol{X}_{\text{d}} -
(\boldsymbol{X}_{\text{d}}:\boldsymbol{R}_{\text{eff}} ) \boldsymbol{R}_{\text{eff}} }{\|\boldsymbol{X}_{\text{d}} -
(\boldsymbol{X}_{\text{d}}:\boldsymbol{R}_{\text{eff}} ) \boldsymbol{R}_{\text{eff}} \| }.
\end{equation}
%Thus, the gradient of the angle $\theta$ has the following properties: its norm
%equals $1 / \|\boldsymbol{\sigma}^{\text{D}}_{\text{eff}}\|$,
%it is orthogonal to the radial direction $\boldsymbol{R}_{\text{eff}}$ and it
%is pointing away from $\boldsymbol{X}_{\text{d}}$ .
Further, differentiating \eqref{overDe3} we obtain
\begin{equation}\label{App5}
\frac{ \partial f(\boldsymbol{\sigma}, \boldsymbol{X}_{\text{k}}, \boldsymbol{X}_{\text{d}}, R) }
{\partial \boldsymbol{\sigma}} = \sqrt{\frac{2}{3}} (K_0 + R)
\frac{\partial \bar{f}\big(\vec{y}(\| \boldsymbol{\sigma}^{\text{D}}_{\text{eff}} \|, \theta ),\alpha \big) }
{\partial \boldsymbol{\sigma}},
\end{equation}
where
\begin{equation}\label{App6}
\vec{y}(\| \boldsymbol{\sigma}^{\text{D}}_{\text{eff}} \|, \theta ) =
\frac{\sqrt{3/2} \ \| \boldsymbol{\sigma}^{\text{D}}_{\text{eff}}  \|}{K_0+R}
(\cos(\theta), \sin(\theta)).
\end{equation}
It follows from \eqref{App6} that
\begin{equation}\label{App7}
\frac{\partial \vec{y} }
{\partial \| \boldsymbol{\sigma}^{\text{D}}_{\text{eff}}  \|} =
\frac{ \sqrt{3/2}}{K_0+R}
(\cos(\theta), \sin(\theta)), \
\frac{\partial \vec{y} }
{\partial \theta} = \frac{ \sqrt{3/2} \ \| \boldsymbol{\sigma}^{\text{D}}_{\text{eff}}  \|}{K_0+R}
(-\sin(\theta), \cos(\theta)).
\end{equation}
Next, note that the degree of distortion $\alpha$ is to be held constant as well, since it is a unique function
of $\boldsymbol{X}_{\text{d}}$.
Thus, using the chain rule we get from \eqref{App5}
\begin{multline*}
\frac{ \partial f(\boldsymbol{\sigma}, \boldsymbol{X}_{\text{k}}, \boldsymbol{X}_{\text{d}}, R) }
{\partial \boldsymbol{\sigma}} = \\
 \sqrt{\frac{2}{3}} (K_0 + R) \Bigg[ \Bigg(
\frac{\partial \bar{f}\big(\vec{y},\alpha \big) }
{\partial \vec{y}} \cdot
\frac{\partial \vec{y} }
{\partial \| \boldsymbol{\sigma}^{\text{D}}_{\text{eff}}  \|} \Bigg)
\frac{ \partial  \| \boldsymbol{\sigma}^{\text{D}}_{\text{eff}} \|}
{\partial \boldsymbol{\sigma}} +
\Bigg( \frac{\partial \bar{f}\big(\vec{y},\alpha \big) }
{\partial \vec{y}} \cdot
\frac{\partial \vec{y} }
{\partial \theta} \Bigg)
\frac{ \partial  \theta}
{\partial \boldsymbol{\sigma}} \Bigg].
\end{multline*}
Substituting \eqref{App1}, \eqref{App4}, and \eqref{App7} into this result we get
for $\theta \neq 0$
\begin{multline}\label{App9}
\frac{ \partial f(\boldsymbol{\sigma}, \boldsymbol{X}_{\text{k}}, \boldsymbol{X}_{\text{d}}, R) }
{\partial \boldsymbol{\sigma}} =
 \Big( \frac{\partial \bar{f}\big(\vec{y},\alpha \big) }
{\partial \vec{y}} \cdot
(\cos(\theta), \sin(\theta)) \Big) \boldsymbol{R}_{\text{eff}} + \\
\Big( \frac{\partial \bar{f}\big(\vec{y},\alpha \big) }
{\partial \vec{y}} \cdot
(\sin(\theta), -\cos(\theta)) \Big)
\frac{\boldsymbol{X}_{\text{d}} -
(\boldsymbol{X}_{\text{d}}:\boldsymbol{R}_{\text{eff}} ) \boldsymbol{R}_{\text{eff}} }{\|\boldsymbol{X}_{\text{d}} -
(\boldsymbol{X}_{\text{d}}:\boldsymbol{R}_{\text{eff}} ) \boldsymbol{R}_{\text{eff}} \| }.
\end{multline}
Here, the gradient of the non-dimensional overstress is computed by \eqref{ApA8}.
Finally, the normality vector tends to the radial direction as $\theta \rightarrow 0$. Therefore
\begin{equation*}
\frac{ \partial f(\boldsymbol{\sigma}, \boldsymbol{X}_{\text{k}}, \boldsymbol{X}_{\text{d}}, R) }
{\partial \boldsymbol{\sigma}} = \boldsymbol{R}_{\text{eff}} \quad \text{for} \ \theta = 0.
\end{equation*}

%% References with bibTeX database:

\bibliographystyle{elsarticle-harv}
%\bibliography{<your-bib-database>}

%% Authors are advised to submit their bibtex database files. They are
%% requested to list a bibtex style file in the manuscript if they do
%% not want to use elsarticle-harv.bst.

%% References without bibTeX database:

% \begin{thebibliography}{00}

%% \bibitem must have one of the following forms:
%%   \bibitem[Jones et al.(1990)]{key}...
%%   \bibitem[Jones et al.(1990)Jones, Baker, and Williams]{key}...
%%   \bibitem[Jones et al., 1990]{key}...
%%   \bibitem[\protect\citeauthoryear{Jones, Baker, and Williams}{Jones
%%       et al.}{1990}]{key}...
%%   \bibitem[\protect\citeauthoryear{Jones et al.}{1990}]{key}...
%%   \bibitem[\protect\astroncite{Jones et al.}{1990}]{key}...
%%   \bibitem[\protect\citename{Jones et al., }1990]{key}...
%%   \harvarditem[Jones et al.]{Jones, Baker, and Williams}{1990}{key}...
%%

% \bibitem[ ()]{}

% \end{thebibliography}

\end{document}